\documentclass[EJP,preprint]{ejpecp} 

\usepackage[utf8]{inputenc}

\usepackage{graphicx}
\usepackage{enumitem}
\usepackage[english]{babel}
\usepackage{bm}
\usepackage{pdfpages}
\usepackage{float} 
\usepackage{latexsym}
\usepackage{verbatim} 
\usepackage{color}
\usepackage{soul}
\usepackage{bbm}
\usepackage{relsize}
\usepackage{cancel}
\usepackage{units}
\usepackage{marginnote}
\usepackage[
]{hyperref}
\usepackage{float}
\usepackage{longtable}

\SHORTTITLE{Weak convergence for the Kingman coalescent}

\TITLE{Weak convergence of the scaled jump chain and  number of mutations of the Kingman coalescent}

\AUTHORS{%
  Martina~Favero\footnote{University of Warwick, United Kingdom.
    \EMAIL{martina.favero@warwick.ac.uk}}
  \and 
  Henrik~Hult\footnote{KTH, Royal Institute of Technology, Sweden. \EMAIL{hult@kth.se} }}

\AMSSUBJ{60J90} 
\AMSSUBJSECONDARY{60F05 ; 92D15}

\SUBMITTED{September 14, 2023} 
\ACCEPTED{April 17, 2024} 

\VOLUME{0}
\YEAR{2024}
\PAPERNUM{0}
\DOI{10.1214/YY-TN}

\ABSTRACT{
The Kingman coalescent  is a fundamental  process in population genetics modelling the ancestry of a  sample of individuals backwards in time. 
In this paper, in a large-sample-size regime,
we study asymptotic properties of the coalescent  under neutrality and a general finite-alleles mutation scheme, i.e. including both parent independent and parent dependent mutation. 
In particular, we consider a  sequence of Markov chains that is related to the coalescent and consists of block-counting and mutation-counting components. We show that these components, suitably scaled,  converge weakly to 
deterministic components and Poisson processes with varying intensities, respectively. 
Along the way,  we develop a novel approach, based on a change of measure, to generalise the convergence result from the parent independent to the parent dependent mutation setting, 
in which several crucial quantities are not known explicitly. 
}

\newcommand{\E}[1]{\mathbb{E}\left[#1\right]}
\newcommand{\Prob}[1]{\mathbb{P} \left( #1\right) }
\newcommand{\x}{\textbf{y}}
\newcommand{\y}{\textbf{y}}
\newcommand{\z}{\textbf{z}}
\newcommand{\vv}{\textbf{v}}

\newcommand{\hist}{\boldsymbol{\mathcal{H}}}
\newcommand{\kc}{\textbf{H}}
\newcommand{\e}{\textbf{e}}

\newcommand{\Y}{\textbf{Y}}

\newcommand{\M}{\textbf{M}}
\newcommand{\Z}{\textbf{Z}}
\newcommand{\m}{\textbf{m}}
\newcommand{\w}{\textbf{w}}
\newcommand{\PP}{\textbf{P}}
\newcommand{\Q}{\textbf{Q}}

\newcommand{\metric}{\psi}

\newcommand{\norm}[1]{\left\lVert#1\right\rVert_1}
\newcommand{\norma}[1]{\left\lVert#1\right\rVert_2}
\newcommand{\norminfty}[1]{\left\lVert#1\right\rVert_\infty}

\newcommand{\nth}{^{(n)}}
\newcommand{\tri}{^{\Delta}}
\newcommand{\zero}{\textbf{0}}

\definecolor{airforceblue}{rgb}{0.36, 0.54, 0.66}

\definecolor{aliceblue}{rgb}{0.94, 0.97, 1.0}

\definecolor{oldlace}{rgb}{0.99, 0.96, 0.9}
\definecolor{burntorange}{rgb}{0.8, 0.33, 0.0}

\definecolor{honeydew}{rgb}{0.94, 1.0, 0.94}
\definecolor{forestgreen}{rgb}{0.13, 0.55, 0.13}

\begin{document}

\section{Introduction}

The Kingman coalescent \cite{kingman1982b} is a classical stochastic process that models the genealogical tree of a sample of individuals. 
With the aim of including  various genetic forces, generalisations of this pivotal model resulted in the establishment of other well-known models, such as the ancestral selection graph \cite{Krone1997,Neuhauser1997}, the ancestral recombination graph \cite{griffiths1991,griffiths1997}, the $\Lambda-$coalescent \cite{pitman1999,sagitov1999}, 
and the $\Xi-$coalescent \cite{mohle2001, schweinsberg2000}.
Coalescent theory is also closely related to urn models. In fact, a P\'olya-like urn structure can be embedded in the coalescent process by matching  lineages and  balls \cite{wakeley2020}.

In this paper, we consider the Kingman  coalescent with a finite-allele mutation scheme 
evolving under neutrality (i.e. no type  has a selective advantage), 
and study its asymptotic behaviour as the size of the sample grows to infinity. 
Beside the purely mathematical significance of deriving new properties of a classical model, our analysis is inspired by applied statistical problems which have recently drawn attention to large-sample-size regimes in population genetics and are briefly mentioned in the following.

Coalescent models are often used for inference, combined with Monte Carlo methods to approximate the likelihood of an observed genes configuration.
In particular,  
importance sampling algorithms \cite{griffiths1994simulating,stephens2000,deiorio2004}
and 
sequential Monte Carlo methods \cite{Koskela2018}
have been developed 
based on simulating backwards the genealogy of the sample (extensions based on generalisations of the coalescent exists \cite{stephens2003, birkner2008, griffiths2008, hobolth2008, birkner2011, Koskela2015}).
However, it is known empirically that these methods do not scale well  with increasing sample size \cite{kelleher2016}. 
Another issue is that  
the coalescent approximation is only accurate when the sample size is sufficiently smaller than the effective population size. Otherwise, as explained in \cite{bhaskar2014}, substantial differences between the  original models and the coalescent approximation could arise. 
Studying  large-sample-size asymptotic properties of the coalescent 
provides a tool both for the theoretical asymptotic analysis of inference methods, which could guide their improvement,  and  for the analysis of the appropriateness of the coalescent as sample sizes in modern studies in genetics grow rapidly.  

The  contribution of the paper is twofold.
The  main  result (Theorem \ref{thm:weak_conv}) is the convergence  of a sequence of Markov chains related to the coalescent, which consists of  block-counting and  mutation-counting components and is described in Section \ref{sect:outline}. We show that these components, suitably scaled, converge to deterministic components and Poisson processes with varying intensities, respectively.
The convergence result is first proved under the assumption of  parent independent mutation (PIM), i.e. the type of the mutated offspring does not depend on the type of the parent. 
We then develop a novel approach, based on a change of measure and presented in Section \ref{sect:change_measure}, to remove the PIM assumption and generalise the convergence result. The change-of-measure approach (Theorem \ref{thm:change_measure}) constitutes the other main result of the paper. In fact, to the best of our knowledge, it provides the first tool in population genetics to transfer results that are obtained under the widely-spread PIM assumption to the general mutation setting. This setting  is notoriously more challenging and often neglected in the literature due to several crucial quantities not being known explicitly.
Other technical challenges for the convergence proofs are addressed in Section 4 by constructing a technical framework that  allows employing the classical tools developed by  Ethier and Kurtz  \cite{Ethier1986}.
Section \ref{sect:conv} contains the convergence proofs. 

\section{Outline and main result}
\label{sect:outline}

The  coalescent describes the genealogy of a sample of individuals by describing the evolution of their ancestral lineages which coalesce and mutate.
Time evolves backwards from time $0$, at which the sample is taken, until only one lineage is left, i.e. the most recent common ancestor of the sample is reached. 
Here, lineages are \textit{typed}, that is, types are 
assigned to lineages, rather than superimposed afterwards, as the coalescent evolves backwards conditioned on the initial sample, as in e.g.  \cite{griffiths1994simulating,Etheridge2009,favero2020, favero2021}. We refer to this version of the coalescent as the typed conditional coalescent. 

We consider the sample, taken from a population at stationarity, to be of  the form $n\y_0\nth$, with $\y_0\nth \in \frac{1}{n}\mathbb{N}^d\setminus\{\boldsymbol{0}\} $, where components represent the numbers of individuals of the various possible  types $1,\dots,d$. 
While the typed conditional coalescent models the genealogy of the sample, as described below in details, the frequencies of types in the population are modelled by the Wright-Fisher diffusion. The sample thus corresponds to a multinomial draw from the stationary distribution of the Wright Fisher diffusion. This close connection between the typed conditional coalescent and the Wright-Fisher diffusion finds its explanation in their duality relationship, see e.g.   \cite{barbour2000,Etheridge2009, jenkins2016, favero2021,favero2023}.

Besides the number of types $d$, the key parameters for both the coalescent and the diffusion are the mutation rate $\theta$ and the mutation probability matrix $P=(P_{ij})_{i,j=1}^d$, with $P_{ij}$ being the probability that the type changes from $i$ to $j$ when a mutation occurs. Throughout the paper, as it is customary, $P$ is assumed to be irreducible,  to ensure the existence of a unique stationary distribution. 

\begin{figure}[]
    \centering
    \includegraphics[width=0.4\textwidth]{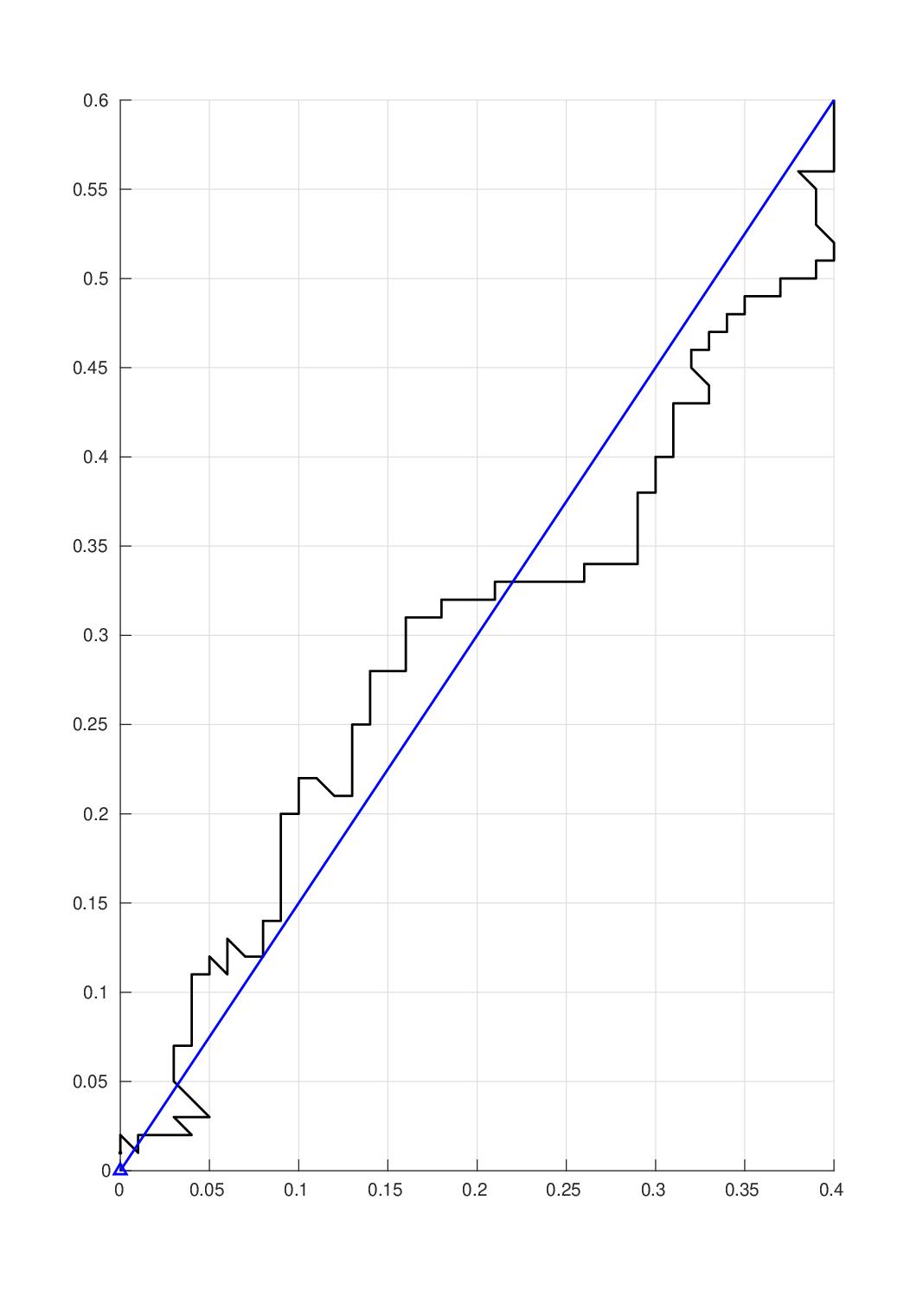}
    \includegraphics[width=0.4\textwidth]{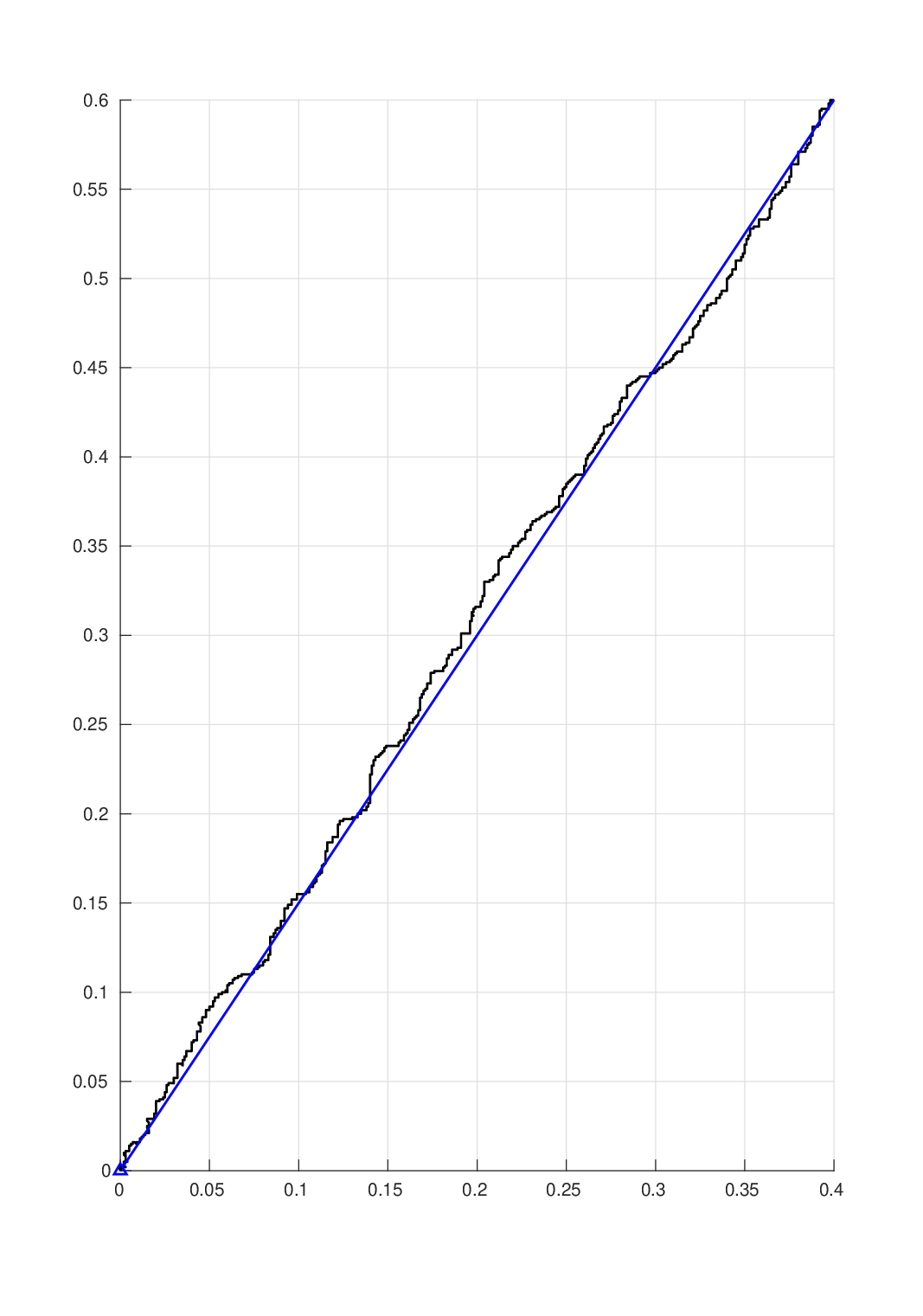}
    \caption{A realisation of the Markov chains $\Y^{(100)}$ (left) and $\Y^{(1000)}$ (right) with starting point $\y_0=(0.4,0.6)$. The limiting process $\Y$ is represented by the blue line. Parameters: mutation rate $\theta=4$, mutation probabilities $Q_1=Q_2=0.5$, number of types $d=2$.
    }
    \label{fig:limit}
    \end{figure}

Instead of considering the whole tree structure formed by the coalescing lineages, we count its \emph{blocks}, that is, the  numbers of lineages of each type over time. 
Furthermore, time is discretised: the original coalescent evolves in continuous time, as well as its block-counting process,  we consider instead the corresponding block-counting jump chain, which jumps to the same states in discrete steps. The $k^{th}$ state of the discrete-time chain corresponds to the state of the continuous-time process right after its $k^{th}$ jump ($k^{th}$ coalescence or mutation event).
The block-counting jump chain of the conditional typed coalescent is denoted by $\kc\nth=\{\kc\nth(k)\}_{k\in \mathbb{N}}$, given the initial sample  $\kc\nth(0)=n \y_0\nth$ and evolving backwards in time until the most recent common ancestor is reached. Thus, $H\nth_i(k)$ denotes the  number of lineages of type $i$ after $k$ jumps in the ancestral history.  
The transitions of the Markov chain $\kc\nth$ are precisely defined in the following, after a scaling.
Let $\Y\nth=\{\Y\nth(k)\}_{k\in \mathbb{N}}\subset\frac{1}{n}\mathbb{N}^d\setminus\{\boldsymbol{0}\}$ be defined as 
    \begin{align*}
    \Y\nth (k)=\frac{1}{n}\kc\nth(k)=\frac{1}{n}\left(H_1(k), \dots, H_d(k)\right),    
    \end{align*}  
see Figure \ref{fig:limit}. 
Given a current state $\Y\nth(k)=\y\in \frac{1}{n}\mathbb{N}^d\setminus\{\boldsymbol{0}\} $,  the next state $\Y\nth(k+1)$ is of the form $\y-\frac{1}{n}\vv$, where,   denoting the $d$-dimensional unit vector along the $j^{th}$ dimension by  $\e_j$, 
    \begin{itemize}
        \item $\vv=\e_i, i=1,\dots,d$, with probability $\rho\nth(\e_i \mid \y)$, which corresponds to the coalescence of two lineages of type $i$;
        \item $\vv=\e_i-\e_j, i,j=1,\dots,d$, with probability $\rho\nth(\e_i-\e_j \mid \y)$, which corresponds to the mutation of a lineage of type $j$ to type $i$ forwards in time (or, equivalently,  $i$ to  $j$ backwards in time).
    \end{itemize}
The backwards transition probabilities $\rho\nth$ are generally not known in an explicit form, unless mutations are parent independent, 
see e.g. \cite{griffiths1994simulating,stephens2000,deiorio2004,griffiths2008} for a more detailed description. Implicit and explicit expressions for $\rho\nth$, in the general and PIM case respectively, are provided in  \eqref{eq:rho} in  Appendix \ref{appendix:transition_prob} and in  \eqref{eq:rho_PIM} in Section \ref{sect:conv_generators}.

Furthermore, along with $\Y\nth$, consider the process $\M\nth=\{\M\nth(k)\}_{k\in \mathbb{N}}\subset \mathbb{N}^{d^2}$ which counts the number of mutations of each type occurring in the genealogy. That is, $\M\nth (k)=(M\nth_{ij}(k))_{i,j=1}^d$, where $M\nth_{ij}(k)$ is the cumulative number of mutations from type $i$ to type $j$ (forwards, or $j$ to $i$ backwards) that have occurred in the genealogy during the first $k$ steps, i.e. 
    \begin{align*}
    M_{ij}\nth (k)= \sum_{k'=0}^{k-1}
    \mathbb{I}_{\{\Y\nth(k')-\Y\nth(k'+1)=\frac{1}{n}(\e_j-\e_i)\}},
    \end{align*} 
and $M_{ij}(0)=0$.

In this paper, we study the asymptotic behaviour of the sequence 
$ \Z\nth=  (\Y\nth ,\M\nth )$, as $n\to\infty$. 
Note that, the sequence $\M\nth$ is not only interesting for the purpose of studying an additional asymptotic property of the coalescent, but it is also crucial because it naturally appears in our change-of-measure argument, as explained in the next sections. 
The main convergence result is stated in the following theorem.

\begin{theorem}
\label{thm:weak_conv}    
Let $\y_0\in \mathbb{R}_+^d$ and $\y_0\nth\in\frac{1}{n}\mathbb{N}^d\setminus\{\boldsymbol{0}\}, n\in \mathbb{N} $. 
Assume $\Y\nth(0)= \y_0\nth$, and 
    $
    \y_0\nth \to \y_0    
    $
as $n\to\infty$. Then,  for all $t\in [0,\norm{\y_0})$, as $n\to\infty,$ the sequence of processes 
    $ \tilde{\Z}\nth= \{ \Z\nth(\lfloor{sn}\rfloor{} )  \}_{s\in[0,t]}$
converges weakly to the process 
    $ 
    \Z= \{ (\Y(s ) ,\M(s)) \}_{s\in[0,t]} \subset \mathbb{R}_+^d \times \mathbb{N}^{d^2}
    $, 
where $\Y$ is the deterministic process defined by
    \begin{align*}
    \Y(s)= \frac{\y_0}{\norm{\y_0}} \left(\norm{\y_0} -s \right),
    \end{align*}
and $\M=(M_{ij})_{i,j=1}^d$ , with $M_{ij}=\{M_{ij}(s)\}_{s\in[0,t]}$ being independent time-inhomogeneous Poisson processes with intensities  
    \begin{align*}
     \lambda_{ij}(\Y(s))= \frac{\theta P_{ij} Y_i(s)}{\norm{\Y(s)}^2}.
    \end{align*}
\end{theorem}
Here converging weakly means converging  in the Skorokhod space  $D_{\mathbb{R}_+^d \times \mathbb{N}^{d^2}} [0,t]$. That is, for  any 
bounded continuous real-valued function $g$ on  $D_{\mathbb{R}_+^d \times \mathbb{N}^{d^2}} [0,t]$, it yields
    \begin{equation*}
    \lim_{n\to\infty} 
    \E{ g\left(
    \{ \tilde{\Z}\nth(s )  \}_{s\in[0,t]}
    \right)}
    = \E{g\left(
    \{ \Z(s )  \}_{s\in[0,t]}
    \right)}.
    \end{equation*}
We now provide an intuitive explanation of Theorem \ref{thm:weak_conv} and illustrate the two main challenges that need to be addressed to obtain a rigorous proof. 
Define the operator $A\nth$,  associated to the Markov chain $\Z\nth$, as
    \begin{align}
    \nonumber
    A\nth f(\x,\m) 
    &=
    n \E{ f(\Z\nth(k+1))-f(\Z\nth(k))\mid \Z\nth(k)=(\y,\m)  }
    \\\nonumber  
    &=
    \sum_{j=1}^{d} 
    n\left[f\left(\x -\frac{1}{n}\e_j,\m\right) -f(\x,\m)\right]
    \rho\nth (\e_j \mid \x) 
    \\
    \label{eq:An}
    & \quad +
    \sum_{i,j=1}^{d}
    \left[f\left(\x -
    \frac{1}{n}\e_j+\frac{1}{n}\e_i,\m+\e_{ij}\right) -f(\x,\m)\right]
    n \rho\nth (\e_j-\e_i \mid \x),
    \end{align}
where $\e_{ij}$ is the unit matrix having $1$ in position $ij$ and $0$ everywhere else. 
Recent results   on the asymptotic behaviour of $\rho\nth$ \cite{favero2020}, show that, if $\x\nth\to\x\in\mathbb{R}_{+}^d$, as $n\to\infty$, then, for $i,j=1,\dots,d$,
    \begin{align}
    \label{eq:transition_prob_limit}
    \rho \nth (\e_j \mid \x\nth)
    \xrightarrow[n\to \infty]{}
    \frac{y_j}{\norm{\x}},
    \quad\quad 
    n \rho\nth(
    \e_j-\e_i  \mid \x\nth
    )
    \xrightarrow[n\to\infty]{}
    \lambda_{ij}(\x)
    .
    \end{align} 
Therefore, by evaluating \eqref{eq:An} at $\y\nth$,  and letting $\y\nth\to\y\in \mathbb{R}_{+}^d$ as $n\to\infty$, 
 it is straightforward to see that the sequence of operators $A\nth$  converges to, in some sense to be properly defined,  
    \begin{equation}
	\label{eq:A}
	A f(\x,\m)= 
	-\left\langle 
	\nabla_\x f(\x,\m),\frac{\x}{\norm{\x}}
	\right\rangle
	+
	\sum_{i,j=1}^d
	\left[
	f(\x,\m+\e_{ij})-f(\x,\m)
	\right]
	\lambda_{ij}(\y).
	\end{equation}
The operator $A$, defined in \eqref{eq:A} above, can be proven to be, as the intuition suggests, the infinitesimal generator of the limiting process of Theorem \ref{thm:weak_conv}, see Appendix \ref{appendix:generator_semigroup}.
In order to make rigorous the convergence of generators sketched above to then prove convergence of the corresponding processes, we need to define the appropriate space of functions, in the domain of the generator $A$, for which the convergence holds. Two main challenges arise. 

The first challenge is given by the lack of an explicit expression for the transition probabilities $\rho\nth$ in the parent dependent mutation case. 
Even though the pointwise asymptotic behaviour \eqref{eq:transition_prob_limit}  of the transition probabilities is known, this is not enough to prove convergence of generators, 
as evident in Section \ref{sect:conv_generators}, because additional information concerning uniform convergence would be needed.  
We address this issue by first proving the convergence  in the PIM case, for which explicit expressions are available, and by then employing a novel approach based on a change of measure to generalise the result from a PIM setting to a general setting, as explained in details in Section \ref{sect:change_measure}. 
It is important to point out that the counting-mutations sequence  $\M\nth$ would naturally appear through the change of measure, even if only the block-counting sequence $\Y\nth$ was considered to begin with (equation \eqref{eq:likelihood} in Section \ref{sect:change_measure}). This strongly motivates the study of the asymptotic behaviour of $\M\nth$, without which the generalisation to the parent dependent setting would not be possible.

The second challenge, which arises already in the simpler PIM setting, is that, if the components of $\y$ are allowed to be arbitrarily close to zero, the scaled transition probabilities of mutation in \eqref{eq:transition_prob_limit} explode, see  \eqref{eq:rho_PIM} in Section \ref{sect:conv_generators}. 
As more precisely described in Section \ref{sect:explosion}, this issue can be addressed  by constructing an appropriate technical framework which relies on the definition of a suitable metric to avoid components being arbitrarily close to zero and allows employing  Ethier-Kurtz classical tools \cite{Ethier1986}. 

While sections \ref{sect:change_measure} and \ref{sect:explosion} are dedicated, respectively, to the change-of-measure argument  and the technical framework, Section 5 contains the convergence proofs leading to, and including, the proof of Theorem \ref{thm:weak_conv}.

\section{From parent independent to parent dependent mutations: a novel change-of-measure approach}
\label{sect:change_measure}

Suppose that the convergence of Theorem \ref{thm:weak_conv} holds 
when mutations are parent independent, that is when the mutation probability matrix $Q=(Q_i)_{i=1}^d$ has identical rows. 
The goal of this section is to illustrate how to generalise the convergence result to a general setting, with  mutation probability matrix $P=(P_{ij})_{i,j=1}^d$. 

More precisely, 
 using subscripts $Q$ and $P$ to  indicate the underlying mutation probability matrices and
fixing a  bounded continuous function $g:D_{\mathbb{R}_+^d \times \mathbb{N}^{d^2}} [0,t]\to \mathbb{R}$, suppose that we know 
    \begin{equation*}
    \lim_{n\to\infty} 
    \E{ g\left(
    \{ \tilde{\Z}_Q\nth(s )  \}_{s\in[0,t]}
    \right)}
    = \E{g\left(
    \{ \Z_Q(s )  \}_{s\in[0,t]}
    \right)},
    \end{equation*}
and that we want to prove  
\begin{equation*}
    \lim_{n\to\infty} 
    \E{ g\left(
    \{ \tilde{\Z}_P\nth(s )  \}_{s\in[0,t]}
    \right)}
    = \E{g\left(
    \{ \Z_P(s )  \}_{s\in[0,t]}
    \right)}.
    \end{equation*}
Our approach consists of rewriting the expectations with respect to $\tilde{\Z}_\PP\nth$ and $\Z_\PP$ as expectations with respect to $\tilde{\Z}_\Q\nth$ and $\Z_\Q$, respectively. That is, we  apply a change of measure to transform the parent dependent mutation setting in a parent independent mutation setting. 
We present the core idea in the following.

Consider the probability of a realisation $\hist(0),\dots,\hist(K), K\in\mathbb{N},$ of  $\kc_P(0), \dots , \kc_P(K)$, given $\kc_P(0)=\hist(0)$, 
    \begin{align*}
    &\Prob{\kc_P(K)=\hist(K),\dots,\kc_P(1)=\hist(1) \mid \kc_P(0)=\hist(0)}
    \\ 
    & \quad \quad =
    \frac{1}{ p_P(\hist(0)) }
    \Prob{\kc_P(K)=\hist(K),\dots,\kc_P(0)=\hist(0)}
    \\
    &\quad \quad =
   \frac{p_P(\hist(K) )}{ p_P(\hist(0)) }
    \prod_{k=0}^{K-1}
    \Prob{\kc_P(k)=\hist(k)\mid \kc_P(k+1)=\hist(k+1)} ,
    \end{align*}
where $p_P(\hist(k))$ is the probability of observing a sample $\hist(k)$ under the unconditioned coalescent, i.e. the sampling probability, and $\Prob{\kc_P(k)=\hist(k)\mid \kc_P(k+1)=\hist(k+1)}$ is the explicitly-known forward transition probability.
More precisely, $p_P$, which can  be expressed explicitly only when mutations are parent independent, can be defined either through a recursion formula or as a multinomial draw from the stationary distribution of the corresponding Wright-Fisher diffusion,
see \cite{griffiths1994simulating} for details. 
The stationary distribution, which is Dirichlet in the PIM case, is also not known explicitly in the general case but exists as long as the mutation probability matrix is irreducible \cite{shiga1981}.  
The boundary condition is assumed to be $p_P(e_i)= \pi_{P,i}, i=1,\dots,d$, where $\boldsymbol{\pi}_P=(\pi_{P,i})_{i=1}^d$ is the invariant distribution of the mutation probability matrix. Despite the implicit formulation, a large-sample-size asymptotic formula  for the sampling probabilities is available in general \cite{favero2020}, which makes the rewriting in the last display useful for the asymptotic analysis.
The forward transition probabilities instead are always explicitly known  \cite{stephens2000}.
Given a configuration $ \kc_P(k+1)=\hist(k+1)$, 
a lineage is chosen uniformly at random, i.e. of type $i$ with probability $\frac{\mathcal{H}_i(k+1)}{\norm{\hist(k+1)}} $. The chosen lineage is split 
 with probability $\frac{\norm{\hist(k+1)}-1}{\norm{\hist(k+1)}-1+\theta}$, and it is mutated to type $j$  with probability  $\frac{\theta}{\norm{\hist(k+1)}-1+\theta} P_{ij}$.
That is, the next configuration is 
$\kc_P(k)= \hist(k+1) + \vv$, where
    \begin{itemize}
    \item $\vv=\e_i, i=1,\dots,d$, with probability $\frac{\mathcal{H}_i(k+1)}{\norm{\hist(k+1)}} 
    \frac{\norm{\hist(k+1)}-1}{\norm{\hist(k+1)}-1+\theta}$;
    \item $\vv=\e_j-\e_i, i,j=1,\dots,d$, with probability 
    $\frac{\mathcal{H}_i(k+1)}{\norm{\hist(k+1)}} 
    \frac{\theta}{\norm{\hist(k+1)}-1+\theta} P_{ij}$.
    \end{itemize}
Note that this description clarifies the equivalence between the coalescent and the P\'olya-like urn model mentioned in the introduction: according to the probabilities above, typed-balls, which correspond to typed-lineages, are duplicated and replaced, instead of split and  mutated respectively  \cite{wakeley2020}. 

It is now possible to derive the likelihood ratio, the Radon–Nikodym derivative,  of the change of measure from parent dependent mutation with mutation probability matrix $P$ to parent independent mutation with mutation probability matrix $Q$. Assuming, without loss of generality, that $Q$ is positive, it yields
    \begin{align}
    \nonumber    &\frac{\Prob{\kc_P(K)=\hist(K),\dots,\kc_P(1)=\hist(1) \mid \kc_P(0)=\hist(0)}}{\Prob{\kc_Q(K)=\hist(K),\dots,\kc_Q(1)=\hist(1) \mid \kc_Q(0)=\hist(0)}}
    \\  \nonumber  &\quad \quad =
    \frac{p_P(\hist(K) )}{ p_Q(\hist(K))}
    \frac{p_Q(\hist(0) )}{p_P(\hist(0))}
    \prod_{k=0}^{K-1}
    \frac{\Prob{\kc_P(k)=\hist(k)\mid \kc_P(k+1)=\hist(k+1)} }{\Prob{\kc_Q(k)=\hist(k)\mid \kc_Q(k+1)=\hist(k+1)} }
    \\  \label{eq:likelihood}
    &\quad \quad =
    \frac{p_P(\hist(K) )}{ p_Q(\hist(K))}
    \frac{p_Q(\hist(0) )}{p_P(\hist(0))}
    \prod_{i,j=1}^d
    \left(
    \frac{P_{ij}}{Q_j}
    \right)^{\sum_{k=0}^{K-1}
    \mathbb{I}_{\{\hist(k)-\hist(k+1)=\e_j-\e_i\}}} .
    \end{align}
Note that the likelihood ratio above does not depend on the whole history, 
but instead only depends on the initial and last configurations, $\hist(K),\hist(0)$, and on the total number of mutations from type $i$ to type $j$, $\sum_{k=0}^{K-1} \mathbb{I}_{\{\hist(k)-\hist(k+1)=\e_j-\e_i\}}$, for $i,j=1,\dots,d$.
This is the key observation for the change-of-measure approach and it is also the reason why it is crucial for the asymptotic analysis to consider, along with the block-counting sequence,  the counting-mutation sequence  which naturally appears in the likelihood ratio.   

In the next theorem, we tailor the previous  calculations to the change of measure from   $\tilde{\Z}_P\nth$ to $\tilde{\Z}_Q\nth$ and derive  similar calculations for the change of measure from $\Z_P$ to $\Z_Q$. Note that an analogous  result can be obtained for any process that is adapted to the block-counting jump-chain of the Kingman coalescent. 

\begin{theorem}
\label{thm:change_measure}
Let $g:D_{\mathbb{R}_+^d \times \mathbb{N}^{d^2}} [0,t]\to \mathbb{R}$ be a bounded continuous function. 
Assume $\tilde{\Y}_P\nth(0)=\tilde{\Y}_Q\nth(0)=\y_0\nth\in \frac{1}{n}\mathbb{N}^d\setminus\{\boldsymbol{0}\}$ 
.
Then, for all $t\geq0$, 
    \begin{align}
    \label{eq:exp.n.rewriting}
    &\E{g\left(\{\tilde{\Z}\nth_P(s)\}_{s\in [0,t]}\right)}
    =     \E{
    g\left(\{\tilde{\Z}_Q\nth(s)\}_{ s\in [0,t]}\right)
    r\nth_{P,Q}\left(\tilde{\Y}\nth_Q(t),\y\nth_0\right)
    c_{P,Q}\left(\tilde{\M}\nth_Q(t)\right)
    },
    \end{align}
where
    \begin{align*}
    r\nth_{P,Q}(\y\nth,\y_0\nth)=
    \frac{p_P(n\y\nth)}{p_Q(n\y\nth)}\frac{p_Q(n\y_0\nth)}{p_P(n\y\nth_0)},
    \quad \y\nth \in \frac{1}{n}\mathbb{N}^d\setminus\{\boldsymbol{0}\},
    \end{align*}
and 
    \begin{align*}
    c_{P,Q}(\m)=\prod_{i,j=1}^{d}\left(\frac{P_{ij}}{Q_j}\right)^{m_{ij}},
    \quad \m \in \mathbb{N}^{d^2}.
    \end{align*}
Furthermore, for all $t\in[0,\norm{\y_0})$, 
    \begin{equation}
    \label{eq:exp.rewriting}
    \begin{aligned}
    \E{g\left(\{{\Z}_P(s)\}_{  s\in [0,t]}\right)}
    =
    \E{g\left(\{{\Z}_Q(s)\}_{ s\in [0,t]}\right)c_{P,Q}\left({\M}_Q(t)\right)	}.
    \end{aligned}
    \end{equation}
\end{theorem}
\begin{proof}
By definition, $\{\tilde{\Z}\nth_P(s)\}_{s\in [0,t]}$ is a function of $n$ and $\kc_P(0),\dots, \kc_P(\lfloor{tn}\rfloor{} ) $. Thus, we can simply use \eqref{eq:likelihood} with $K=\lfloor{tn}\rfloor{} $ to obtain \eqref{eq:exp.n.rewriting}. 

Since the first component of $\Z_P=(\Y,\M_P)$ and of $\Z_Q=(\Y,\M_Q)$ are deterministic and do not depend on the mutation probability matrix, 
it is sufficient to focus on the second components. 
Using that 
    \begin{align*}
    \Prob{\M_P(t)=\m \mid \Y(0)=\y_0}
    =\prod_{i,j=1}^{d} 
    \frac{1}{m_{ij}!}
    \left( \int_0^t \lambda_{P,ij} (\Y(u)) du \right)
    ^{m_{ij}}
    e^{-\int_0^t \lambda_{P,ij} (\Y(u)) du}
    \end{align*}
and noting that  $\int_0^t \lambda_{P,ij} (\Y(u)) du= \frac{\theta P_{ij} y_{0,i}}{\norm{\y_0}} \int_0^t \frac{1}{\norm{\y_0}-u } du $, we obtain
    \begin{align*}
    \frac{\Prob{\M_P(t)=\m \mid \Y(0)=\y_0}}
    {\Prob{\M_Q(t)=\m \mid \Y(0)=\y_0}}
    = 
    c_{P,Q}(\m).
    \end{align*}
In fact, $c_{P,Q}(\M_Q(t))$ can be interpreted as the exponential martingale associated to the change of measure that transforms the distribution of the Poisson processes 
    $
    \{\M_P(s)\}_{ s\in [0,t]}
    $
into the distribution of the Poisson processes 
    $
    \{\M_Q(s)\}_{ s\in [0,t]}
    $,	
see e.g. \cite{kurtz2001}. This proves \eqref{eq:exp.rewriting} and completes the proof.
\end{proof}

\section{Handling the explosion: technical framework}
\label{sect:explosion}

In this section, we construct a  technical framework that allows us to address the problem caused by  
the scaled transition probabilities of mutation in \eqref{eq:transition_prob_limit} exploding near 
the boundary 
    $
    \Omega_0 :=\{\y=(y_1,\dots,y_d) : y_i=0 \text{ for some } i \}
    $.  
In such a framework, we can rigorously state the convergence of generators (Proposition \ref{prop:conv_generators} below)  and  we can show the weak convergence of the sequence $\Z\nth$ to the process $\Z$  by relying on classical results \cite{Ethier1986}, as  explained  briefly below and proved in detail in the next section.

A natural state-space of the limiting process $\Z$  is 
$\mathbb{R}_+^d \times \mathbb{N}^{d^2} $, equipped with the Euclidean metric.
Consider instead $E_1=(0,\infty]^d$, equipped with the metric $\psi_1(\y_1,\y_2):=\norma{\frac{1}{\y_1}-\frac{1}{\y_2}}$, where the inversion is component-wise and the inverse of $\infty$ is by convention $0$;  and consider  
$E= E_1 \times \mathbb{N}^{d^2} $, equipped 
with the product  metric $\metric:=\metric_1\oplus  \norma{\cdot}$. For completeness, we naturally define $\Z$ to be constantly equal to its initial state, if one of the components of $\Y(0)$ is equal to $\infty$, although this type of initial condition does not appear in the relevant setting of Theorem \ref{thm:weak_conv}.

In the new state-space $E_1$, the roles of $0$ and $\infty$ are reversed for each component, whilst, away from $0$ and $\infty$,
the metric $\metric_1$ is equivalent to the Euclidean metric.  In particular, compact sets are bounded away from $\Omega_0$, 
 and thus  functions with compact support are equal to zero near $\Omega_0$. See Appendix \ref{appendix:spaces} for more details on the metric  space $(E,\metric)$.

Let $C_c^{\infty}(E)$ be the space of real-valued continuous functions on $(E,\metric)$ with compact support, equipped with the supremum norm, which is analysed in details in Appendix \ref{appendix:spaces}. 
This is precisely the space of functions that we need for the convergence of generators, which is stated in the next proposition. 

Let $E\nth= \frac{1}{n}\mathbb{N}^d\setminus\{\boldsymbol{0}\}\times \mathbb{N}^{d^2} $ 
be the state space of $\Z\nth$ 
and let 
$\eta_n $ map any function on $E$ into its restriction on $E\nth$, with value zero on $\Omega_0 \times \mathbb{N}^{d^2}$.

\begin{proposition}
\label{prop:conv_generators}
Assume parent independent mutations.
Let $A\nth $ be the operator,  defined in
\eqref{eq:An}, associated to the Markov chain $\Z\nth$, and let and $A  $ be the infinitesimal generator, defined in \eqref{eq:A} , associated to the Markov process $\Z$.   Then, for all $f\in C_c^{\infty}(E)$,
    \begin{equation*}
    \label{}
    \lim_{n\to\infty}
    \sup_{(\x,\m)\in E\nth}
    \left|A\nth \eta_n f(\x,\m) - \eta_n A f (\x,\m)\right|
    =0.
    \end{equation*} 
\end{proposition}
Proposition \ref{prop:conv_generators}, which is proved in Subsection \ref{sect:conv_generators}, allows to derive the convergence of the corresponding semigroups, which is proved in Section \ref{sect:conv_semigroups}.
In order to then prove weak convergence (Theorem \ref{thm:weak_conv}), another technical issue needs to be addressed. In fact, as time approaches  $\norm{\y_0}$, $\Y$ approaches the origin, and the intensities of the Poisson processes $\M$ explode, that is, the process $\Z$ exits its  state space $E$ in a finite time. 
This is a classical problem, and the  technical framework introduced in this section is constructed precisely to match the framework of  Ethier and Kurtz \cite{Ethier1986} and to employ their classical results  which allow us to address the issue and finally prove weak convergence.

Following \cite{Ethier1986}, a new point $\Delta$, called the \textit{cemetery point}, is introduced in the state space. 
Let $E^\Delta = E \cup \{\Delta\}$ be the one-point compactification of $E$,  which is itself a metric space, see Appendix \ref{appendix:spaces} and \cite{mandelkern1989} for more details,  and let 
 $\Z^\Delta$ be the extension of $\Z$, corresponding   to $\Z$ before time $\norm{\y_0}$, and being equal to $\Delta$ from time $\norm{\y_0}$ on. 
Furthermore, let $\xi\tri_n$ be the inclusion  map from $E\nth$ into $E\tri$, sending $\Omega_0\times \mathbb{N}^{d^2} $ to $\Delta$. 
Instead of proving directly Theorem \ref{thm:weak_conv}, i.e. 
    $$
    \{\tilde{\Z}\nth(s )\}_{s\in[0,t]} 
    \xrightarrow[n\to\infty]{} 
    \{\Z(s )\}_{s\in[0,t]} 
    \quad \text{in }  D_{S} [0,t],
    $$
for a fixed $t\in[0,\norm{\y_0})$, we prove an extended convergence 
    \begin{align}
    \label{eq:extended_weak_conv}
    \xi_n\tri \left(\{\tilde{\Z}\nth(s )\}_{s \in [0,\infty) }\right)
    \xrightarrow[n\to\infty]{}
    \{\Z\tri(s )\}_{s\in [0,\infty) }
    \quad \text{in }  D_{E\tri} [0,\infty)
    \end{align}
which then implies Theorem \ref{thm:weak_conv}, as shown in Subsection \ref{sect:weak_conv_PIM}.

\section{Convergence proofs}
\label{sect:conv}
In this section, we prove Theorem \ref{thm:weak_conv}, i.e. the weak convergence of the sequence $\{\tilde{\Z}\nth(s )\}_{s\in[0,t]}$ to the process $\{\Z(s )\}_{s\in[0,t]} $. 
The plan for the proof of  is  the following.
First, under the assumption of parent independent mutations,
\begin{itemize}
\item[$\bullet$] [Sect. \ref{sect:conv_generators}]
 the convergence of generators (Proposition \ref{prop:conv_generators}) is proved;
\item[$\bullet$] [Sect. \ref{sect:conv_semigroups}]
 the convergence of the corresponding semigroups is proved;
\item[$\bullet$] [Sect. \ref{sect:weak_conv_PIM}]
 the weak convergence (Theorem \ref{thm:weak_conv}) is proved   by showing the weak convergence  of the extended processes on the extended state space \eqref{eq:extended_weak_conv}, which relies on the technical framework of Section \ref{sect:explosion}.
\end{itemize}
Then, in the general case of possibly parent dependent mutations,
\begin{itemize}
\item[$\bullet$] [Sect. \ref{sect:weak_conv}]
we use the change-of-measure approach (Theorem \ref{thm:change_measure}) to prove the weak convergence (Theorem \ref{thm:weak_conv}).

\end{itemize}


\subsection{Convergence of generators (PIM) }
\label{sect:conv_generators}

In this subsection, we prove Proposition \ref{prop:conv_generators}.
Assume the mutation probability matrix is of the form $Q=(Q_i)_{i=1}^d$. Then, the backwards transition probabilities are explicitly known:
    \begin{align}
    \label{eq:rho_PIM}
    \rho\nth(\vv \mid \x)
    =&
    \begin{cases}
    \frac{y_j(n y_j -1)}{\norm{\y}(n y_j -1 + \theta Q_j)} 
    &\text{  if  }
    \vv=\e_j,
    \quad j=1\dots d,
    \\
    \frac{\theta Q_j y_j (n y_i - \delta_{ij}+ \theta Q_i)}
    {\norm{\y} (n \norm{\y}-1 + \theta ) ( n y_j -1 + \theta Q_j)}
    &\text{  if  }
    \vv=\e_j-\e_i,
    \quad i, j=1\dots d,
    \\
    0 &\text{  otherwise},
    \end{cases}
    \end{align}
see e.g. \cite{deiorio2004} and Appendix \ref{appendix:transition_prob}.
Given $f\in C_c^{\infty}(E)$,  let $$K=\{(\x,\m)\in E : y_j\geq \delta, m_{ij}\leq M, \forall i,j=1,\dots,d\}$$  be the compact set that contains the support of $f$, see Remark \ref{remark:vanishing_compact} in Appendix \ref{appendix:spaces}. Note that $f$ and  all of its  derivatives  are  Lipschitz  continuous with respect to $\metric$  and bounded. From now on we consider  $n>\frac{1}{\delta}$. 

First, we assume that $(\x,\m)\in E\nth \cap K^c$, which implies that $f=Af=0$ in a neighbourhood of $(\x,\m)$, and that the points
    $
    (\x,\m+\e_{ij})
    $
and 
    $
    \left(\x-\frac{1}{n}\e_j ,\m\right)
    $
belong to  $K^c, \forall n\in\mathbb{N}_{>0}, i,j=1,\dots,d$.
If also 
    $
    \left( \x-\frac{1}{n}\e_j+\frac{1}{n}\e_i ,\m+\e_{ij}\right)
    \in K^c, \forall n\in\mathbb{N}_{>0}, i,j=1,\dots,d, 
    $
then, 
$A\nth \eta_n f(\x,\m)=\eta_n A f(\x,\m)=0$. 
Otherwise, it must be that for a unique $i_0$, and for some $n$,  we have that 
$\delta-\frac{1}{n}\leq y_{i_0}< \delta $; while $y_j\geq\delta$, for all $j\neq i_0$, and $m_{ij}\leq M $, for all $i,j=1,\dots,d$. In this case, 
    $$A\nth \eta_n f(\x,\m)
    = \left| \sum_{j\neq i_0} 
    \eta_n f\left(\x -\frac{1}{n}\e_j+\frac{1}{n}\e_{i_0},\m+\e_{i_0j}\right)
    n\rho\nth (\e_j - \e_{i_0} \mid \x) \right|.
    $$
Note that, since $y_j\geq \delta$, for $j\neq i_0$, we have
    $
     n\rho\nth (\e_j - \e_{i_0} \mid \x)
     < 
     \frac{n \theta Q_j}{n \delta -1 + \theta Q_j},
    $
which is bounded by a constant, independent of $\y$. 
Furthermore, 
    $ 
    \left|\eta_n f\left(\x -\frac{1}{n}\e_j+\frac{1}{n}\e_{i_0},\m+\e_{i_0j}\right)\right| 
    $
is bounded by the supremum of $\left|f(\y',\m')\right| $ over those $(\x',\m')\in E\nth$ having $y'_{i}$ such that $\delta-\frac{1}{n}\leq  y'_{i}< \delta$, for some $i$. This supremum does not depend on $(\y,\m)$ and converges to $0$ as $n\to\infty$ since $f$ is continuous with support in $K$.  
  
Assume next that  $(\x,\m)\in E\nth \cap K $, which implies  $\y\in \frac{1}{n}\mathbb{N}^d \setminus\Omega_0$.  It is to be shown that 
the difference in the next display is bounded by some function of $n$, which does not depend on $\x$ and $\m$, and vanishes as $n\to\infty$. The bound is constructed as 
    \begin{align}
    \allowdisplaybreaks\nonumber
    &
    \left| 
    A\nth \eta_n f(\x,\m) 	- \eta_n Af(\x,\m)
    \right| \\ \nonumber
    &\quad = \left| 
    A\nth  f(\x,\m) 	-  Af(\x,\m)
    \right| 
    \\ \nonumber 
    &\quad\leq
    \sum_{j=1}^d
    \left| 
    \frac{f\left(\x- \frac{1}{n}\e_j,\m\right)-f(\x,\m)}{\frac{1}{n}}	
    \rho\nth(\e_j \mid \x) +
    \frac{\partial f}{\partial y_j} (\x,\m)\frac{y_j}{\norm{\x}}
    \right| 
    \\ \nonumber
    &\qquad+
    \sum_{i,j=1}^d
    \Bigg|
    \left[
    f\left(\x- \frac{1}{n}\e_j+ \frac{1}{n}\e_i,\m+\e_{ij}\right)-f(\x,\m)
    \right]
    n \rho\nth(\e_j-\e_i\mid \x) \\ \nonumber
    &\qquad \quad \quad \quad - 
    \left[
    f(\x,\m+\e_{ij})-f(\x,\m)
    \right]
    \lambda_{Q,ij}(\y)
    \Bigg| 
    \\ \nonumber 
    &\quad\leq
    \sum_{j=1}^d
    \left| 
    \frac{f\left(\x- \frac{1}{n}\e_j,\m\right)-f(\x,\m)}{\frac{1}{n}}	
    +\frac{\partial f}{\partial y_j} (\x,\m) 
    \right| 
    \rho\nth(\e_j \mid \x) 
    \\ \nonumber
    &\qquad+
    \sum_{j=1}^d
    \left| 
    \frac{\partial f}{\partial y_j} (\x,\m)
    \right| 
    \left| 
    \frac{y_j}{\norm{\x}} - \rho\nth(\e_j \mid \x) 
    \right| 
    \\ \nonumber
    &\qquad+
    \sum_{i,j=1}^d
    \left|
    f\left(\x- \frac{1}{n}\e_j+ \frac{1}{n}\e_i,\m+\e_{ij}\right)-f(\x,\m+\e_{ij})
    \right|
    n \rho\nth(\e_j-\e_i \mid \x)
    \\ 
    \label{boundAn-A}
    &\qquad+
    \sum_{i,j=1}^d
    \left|
    f(\x,\m+\e_{ij})-f(\x,\m)
    \right|
    \left|
    n \rho\nth(\e_j-\e_i \mid \x)
    -
    \lambda_{ij}(\y)
    \right|. 
    \end{align}
A bound for each term of the sums in the right hand side of (\ref{boundAn-A}) can be obtained as follows.
Applying the mean value theorem yields
    \begin{align*}
    f\left(\x- \frac{1}{n}\e_j,\m\right)-f(\x,\m)=
    -\frac{1}{n}
    \frac{\partial f }{\partial y_j} \left(a\x+ (1-a)\left(\x-\frac{1}{n}\e_j\right),\m\right),
    \end{align*}
for some $a\in(0,1)$,
therefore, using also that 	
    $
    \rho\nth (\e_j \mid \x)\leq 1
    $, 
and that $\frac{\partial f }{\partial y_j} $ is Lipschitz continuous,
    \begin{align} \nonumber
    \allowdisplaybreaks
    &
    \left| 
    \frac{f\left(\x- \frac{1}{n}\e_j,\m\right)-f(\x,\m)}{\frac{1}{n}}	
    +\frac{\partial f}{\partial y_j} (\x,\m)
    \right| 
    \rho\nth (\e_j \mid \x)
    \\ \nonumber&=  
    \left|
    \frac{\partial f}{\partial y_j} \left(a\x+ (1-a)\left(\x-\frac{1}{n}\e_j\right),\m\right)
    -
    \frac{\partial f }{\partial y_j} (\x,\m)
    \right|
    \rho\nth (\e_j \mid \x)  
    \\ \nonumber &\leq
    \norminfty{\nabla^2_\x f} 
    \metric_1\left(  a\x+ (1-a)\left(\x-\frac{1}{n}\e_j\right),\x \right)
    \\ \nonumber &=
    \norminfty{\nabla^2_\x f}
    \frac{1-a}{y_j(n y_j -1+a)}
    \\ \nonumber &\leq
    \norminfty{\nabla^2_\x f}
    \frac{1-a}{\delta(n \delta -1+a)} .
    \end{align}
For the second term in (\ref{boundAn-A}), by the explicit expression of $\rho\nth(\e_j \mid \y)$ in (\ref{eq:rho_PIM}), it follows that,
    \begin{align*}
    \allowdisplaybreaks
    \left| 
    \frac{\partial f}{\partial y_j} (\x,\m)
    \right| 
    \left| 
    \frac{y_j}{\norm{\x}} - \rho\nth(\e_j \mid \x) 
    \right|& =
    \left| 
    \frac{\partial f}{\partial y_j} (\x,\m)
    \right| 
    \frac{y_j}{\norm{\x}}
    \frac{\theta Q_j}{n y_j-1 +\theta Q_j}
    \\
    &\leq
    \norminfty{\frac{\partial f}{\partial y_j} } 
    \frac{\theta Q_j}{n \delta -1 +\theta Q_j}.
    \end{align*}		
For the third term in (\ref{boundAn-A}), since	
$
\rho\nth(\e_j- \e_i \mid \y) \leq \frac{\theta Q_j}{n \delta-1+\theta Q_j}<1
$, and $f$ is Lipschitz continuous,
    \begin{align*}
    \allowdisplaybreaks
    & \left|
    f(\x- \frac{1}{n}\e_j+ \frac{1}{n}\e_i,\m+\e_{ij})-f(\x,\m+\e_{ij})
    \right|
    n \rho\nth(\e_j-\e_i \mid \x)
    \\ &\quad\leq
    \norminfty{\nabla_\x f}
    \metric_1\left(
    \x- \frac{1}{n}\e_j+ \frac{1}{n}\e_i,\x
    \right)
    \\ & \quad 
    = \norminfty{\nabla_\x f} \norma{
    \frac{\e_j-\e_i}{\y(n\y-\e_j+\e_i))}
    }
    \\ & \quad 
    \leq \norminfty{\nabla_\x f}
    \left(   \frac{1}{y_i(n y_i +1)} + \frac{1}{y_j(ny_j-1)}   \right)
    \\ & \quad 
    \leq \norminfty{\nabla_\x f}
    \left(   \frac{1}{\delta(n \delta +1)} + \frac{1}{\delta(n\delta-1)}   \right).
    \end{align*}
Finally, for the last term in (\ref{boundAn-A}), recalling that in the PIM case $\lambda_{ij}(\y)=\frac{\theta Q_j y_i}{\norm{\y}^2}$, 
    \begin{align*}
    \allowdisplaybreaks
    &
    \left|
    f(\x,\m+\e_{ij})-f(\x,\m)
    \right|
    \left|
    n \rho\nth(\e_j-\e_i \mid \x)
    -
    \lambda_{ij}(\y)
    \right|  \\
    &\quad\leq
    2\norminfty{f}
    \frac{\theta Q_j }{\norm{\x}}
    \left|
    \frac{ny_j}{n y_j -1 +\theta Q_j}\frac{n y_i - \delta_{ij}+Q_i}{n \norm{\x}-1+\theta}
    -
    \frac{y_i}{\norm{\x}}
    \right|
    \\
    &\quad\leq
    2\norminfty{f}
    \frac{\text{const}_\theta}{\delta (n\delta-1)}.
    \end{align*}
Since the bound for \eqref{boundAn-A} does not depend on $(\x,\m)$ and vanishes as $n$ goes to infinity, the proof is complete.   

\subsection{Convergence of semigroups (PIM)}
\label{sect:conv_semigroups}

In this subsection, we prove the convergence of the semigroups  associated to the generators of the previous subsection under the assumption of parent independent mutations. 
Take $f\in \hat{C}(E)$, the space of continuous real-valued function on $(E,\metric)$ vanishing at infinity, see Remark \ref{remark:vanishing_compact} for a characterisation.
Let $T\nth$ be the semigroup associated to the discrete-time Markov chain $\Z\nth$, that is,
    \begin{align}
    \label{eq:Tn}
    T\nth \eta_n f(\x,\m) 
    &=
    \sum_{j=1}^{d} 
    \eta_n f\left(\x -\frac{1}{n}\e_j,\m\right)\rho\nth (\e_j \mid \x) \nonumber 
    \\
    & \quad +
    \sum_{i,j=1}^{d} \eta_n f\left(\x -
    \frac{1}{n}\e_j+\frac{1}{n}\e_i,\m+\e_{ij}\right) \rho\nth (\e_j-\e_i \mid \x), 
    \end{align}
and let $T$ be the semigroup associated to the Markov process $\Z$, described in the following. 
Let $\Omega_\infty:=\{\y=(y_1,\dots,y_d) : y_i=\infty \text{ for some } i \}$.  
If $\x\in E_1\cap \Omega_\infty$,  then $T(t) f(\x,\m)=f(\y,\m)$, $\forall t\geq 0$. If $\x\in E_1\setminus \Omega_\infty$,
then, for $t\geq \norm{\x}$,   $T(t) f(\x,\m)=0$, 
whereas, for $t<\norm{\x} $,
    \begin{align}
    \label{semigroupT}
    T(t) f(\x,\m)
    =&
    \sum_{\w \in \mathbb{N}^{d^2}}
    f\left(\x- \frac{\x}{\norm{\x}}t,\m+\w\right) \gamma_{\w}(t,\x), 
    \end{align}
where,
for $\w=(w_{ij})_{i,j=1,\dots,d} \in \mathbb{N}^{d^2}$,
    \begin{align}
    \label{eq:gamma_prob}
    \gamma_{\w}(t,\x)=&
    \prod_{i,j=1}^{d} 
    \frac{\Lambda_{ij}(t,\y) ^{w_{ij}}}{w_{ij}!}
    e^{-\Lambda_{ij}(t,\y)}.
    \end{align}
with
    \begin{align*}
    \Lambda_{ij}(t,\y)=
    \int_0^t \lambda_{ij} \left(\frac{\y}{\norm{\y}} \left(\norm{\y} -u \right)\right) du
    =
    \frac{\theta Q_j y_i}{\norm{\x}} 
    \log\left(\frac{\norm{\x}}{\norm{\x}-t}\right)
    \end{align*}
We now prove that, for all $f\in\hat{C}(E)$, for all $t\geq 0$,
    \begin{equation}
    \label{TntoT}
    \lim_{n\to\infty}
    \sup_{(\x,\m)\in E\nth}
    \left|(T\nth) ^{\lfloor{tn}\rfloor} 
    \eta_n f(\x,\m) - \eta_n T(t) f (\x,\m)\right|
    =0.
    \end{equation} 
First, note that $T\nth $ is a linear contraction and that, as shown in   Appendix \ref{appendix:generator_semigroup}, $\{T(t)\}_{t\geq 0}$ is a strongly continuous contraction semigroup associated to the generator $A$. 
Furthermore, $C^\infty_c(E)$ is a  core for the generator $A$. 
This holds by Proposition 3.3 in \cite[Ch.1]{Ethier1986}, since  $C^\infty_c(E)$ is dense in $\hat{C}(E)$,   $T(t)$ maps $C^\infty_c(E)$ into $C^\infty_c(E)$ and $A $ is  the generator of a strongly continuous contraction semigroup.
Therefore, by Theorem 6.5 in \cite[Ch.1]{Ethier1986},
the convergence of semigroups (\ref{TntoT}) is equivalent to the convergence of generators of Proposition \ref{prop:conv_generators}, which concludes the proof.

\subsection{Weak convergence (PIM)}
\label{sect:weak_conv_PIM}
In this subsection, we  prove Theorem \ref{thm:weak_conv} under the assumption of parent independent mutations. 
We first prove the  result for the extended processes \eqref{eq:extended_weak_conv}, following \cite{Ethier1986}. 
The semigroup $T$, associated to the Markov process $\Z$,  is extended to a conservative semigroup $T\tri$ as in  e.g. \cite[Ch. 4]{Ethier1986}, which is associated to the extended  Markov process $\Z\tri$, by letting, for $f\in \hat{C}(E\tri)=C(E\tri)$,
    \begin{equation}
    \label{Textension}
    T\tri(t) f= f(\Delta)+T(t) (f-f(\Delta)).
    \end{equation}
In the expression above it is implied that $T(t)$ is applied to the restriction to $E$ of the function $f(\cdot)-f(\Delta)$  and that  $T\tri(t) f $ is a function on $E\tri$, equal to $f(\Delta)+T(t) (f-f(\Delta))$ on $E$ and to $f(\Delta)$ on $\Delta$.
By Lemma 2.3 in \cite[Ch.4]{Ethier1986}, the extended semigroup
$\{T\tri(t)\}_{t\geq 0}$ inherits all the properties of $\{T(t)\}_{t\geq 0}$ and additionally it is conservative. Hence, it is a Feller semigroup.

By defining the bounded linear operator $\eta_n\tri f= f\circ \xi_n\tri$, the convergence of semigroups of the previous subsection can be naturally extended to, for all $f\in C(E\tri)$ and    $t\geq0$,     
    \begin{equation*}
    \lim_{n\to\infty}
    \sup_{(\x,\m)\in E\nth}
    \left|(T\nth) ^{\lfloor{tn}\rfloor} 
    \eta_n\tri f(\x,\m) - \eta_n\tri T\tri(t) f (\x,\m)\right|
    =0.
    \end{equation*}
Therefore, by  Theorem 2.12 in \cite[Ch. 4]{Ethier1986}, 
there exists a Markov process $\Z\tri$, associated to the semigroup $T\tri$, with sample paths in $D_{E\tri}[0,\infty)$, and
    $
    \xi_n(\tilde{\Z}\nth)
    \to
    \Z\tri
    $ 
in  $D_{E\tri}[0,\infty)$.
Since $T\tri$ is the extension of the generator $T$ associated to the process $\Z$, $\Z\tri $ coincides with $\Z$ on $E$.

This easily implies Theorem \ref{thm:weak_conv}, as shown in the following. 
Since $\Y(0)=\y_0$ and $ t < \norm{\y_0}$, the limiting process at time $t$ is not at the cemetery point, that is, $\Z\tri(t)=\Z(t)\in E$.
Letting for example $\delta = \frac{1}{2}\min_{i=1,\dots,d}\{ \tilde{Y}_i(t)\}$, which is positive, the  convergence of the extended processes \eqref{eq:extended_weak_conv} implies that the probability of all  components of $\{\tilde{Y}\nth(s)\}_{0\leq s \leq t}$ being larger than $\delta$ converges to $1$. 
Being the paths up to time $t$ bounded away from $\Omega_0$, 
the  convergence of the extended processes \eqref{eq:extended_weak_conv}, restricted to  the paths up to time $t$, 
is equivalent to the convergence of the non-extended processes  written  in terms of the Euclidean metric.

\subsection{Weak convergence (general mutation) }
\label{sect:weak_conv}

In this subsection, we prove Theorem \ref{thm:weak_conv} for a general mutation matrix $P$, 
knowing from the previous subsection that the theorem holds for a PIM mutation matrix $Q$, which we assume positive, and using the change-of-measure approach, i.e. Theorem \ref{thm:change_measure}.

First, we prove the following lemma on the asymptotic properties of the sequences $r_{P,Q}\nth(\tilde{\Y}_Q\nth(t)$ and 
$c_{P,Q}(\tilde{\M}_Q\nth(t))$  which constitute the Radon–Nikodym derivative of the change of measure from $P$ to $Q$ appearing in Theorem \ref{thm:change_measure}.  We use Theorem \ref{thm:weak_conv} for the sequence $\Z\nth_Q$ and the asymptotic results in \cite{favero2020} for the sampling probabilities. 

\begin{remark}
\label{remark:thm}
For the following proofs it is  useful to recall that, 
for a sequence of random variables that converges in distribution, 
convergence of the $\mathcal{L}^1$-norms occurs if and only if the sequence is uniformly integrable. 
\end{remark}

\begin{lemma} 
\label{lemma:CR}	
Let $\y_0\in \mathbb{R}_+^d$ and $\y_0\nth\in\frac{1}{n}\mathbb{N}^d\setminus\{\boldsymbol{0}\}, n\in \mathbb{N} $. 
Assume $ \Y_Q\nth(0)=\y_0\nth$, and 
    $
    \y_0\nth \to \y_0    
    $
as $n\to\infty$. 
Then,  for all $t\in [0,\norm{\y_0})$, 
\begin{enumerate}[label=(\roman*)]
\item
\label{lemma:CR.exp.n.1}
$
\E{c_{P,Q}(\tilde{\M}_Q\nth(t)) r_{P,Q}\nth(\tilde{\Y}_Q\nth(t), \y_0\nth) 	
}
=1
$;
\item
\label{lemma:CR.exp.1}
$
\E{c_{P,Q}(\M_Q(t))}=
1
$;
\item
\label{lemma:CR.R1}
$
r\nth_{P,Q}(\tilde{\Y}_Q\nth(t), \y_0\nth) \to 1 
$
in probability, as $n\to\infty$;
\item 
\label{lemma:CR.ntoC}
$
c_{P,Q}(\tilde{\M}_Q\nth(t)) r\nth_{P,Q}(\tilde{\Y}_Q\nth(t), \y_0\nth) \to c_{P,Q}(\M_Q(t))$ in distribution, as $n\to\infty$;
\item	
\label{lemma:CR.UI}
$
\{c_{P,Q}(\tilde{\M}_Q\nth(t)) r\nth_{P,Q}(\tilde{\Y}_Q\nth(t), \y_0\nth)
\}_{n\in \mathbb{N}}
$
is a  uniformly integrable sequence.
\end{enumerate}

\begin{proof}
To prove
\ref{lemma:CR.exp.n.1} and \ref{lemma:CR.exp.1},
simply apply Theorem \ref{thm:change_measure} with $g=1$.

To prove \ref{lemma:CR.R1}, first recall that, by 
\cite[Thm 4.3]{favero2020}, 
$
n^{d-1} p_P(n\y\nth) \to \norm{\y}^{1-d}\tilde{p}_P(\y)  ,
$
as $\y\nth\to \y$, where $\tilde{p}_P$ is the stationary density of the Wright-Fisher diffusion that is dual to the Kingman coalescent.
Thus
    \begin{align*}
    \lim_{n\to\infty}
    r\nth_{P,Q}(\y\nth,\y_0\nth)
    =
    \frac{\tilde{p}_P\left(\frac{\y}{\norm{\y}}\right)}
    {\tilde{p}_Q\left(\frac{\y}{\norm{\y}}\right)} 
    \frac{\tilde{p}_Q\left(\frac{\y_0}{\norm{\y_0}}\right)}
    {\tilde{p}_P\left(\frac{\y_0}{\norm{\y_0}}\right)}. 
    \end{align*}
Furthermore,  Theorem \ref{thm:weak_conv}, which is already proven under the PIM assumption, implies that 
$\tilde{\Y}\nth_Q(t)$ 
converges in probability to the constant $\Y(t)$, 
which is equal to $\frac{\y_0}{\norm{\y_0}} \left(\norm{\y_0} -t \right)$. 
By combining the two convergences above and using the characterization of convergence in probability in terms of almost surely converging subsequences, it is straightforward to conclude that $r\nth_{P,Q}(\tilde{\Y}_Q\nth(t), \y_0\nth)$ converges in probability to $1$,
since 
    \begin{align*}
    \frac{\Y(t)}{\norm{\Y(t)}}
    =
    \frac{\frac{\y_0}{\norm{\y_0}} \left(\norm{\y_0} -t \right)}{\norm{\frac{\y_0}{\norm{\y_0}} \left(\norm{\y_0} -t \right)}}
    =\frac{\y_0}{\norm{\y_0}}.
    \end{align*}

Proving  \ref{lemma:CR.ntoC} is rather straightforward. In fact,
since  $c_{P,Q}$ is a continuous function, and since, by Theorem \ref{thm:weak_conv}, $\tilde{\M}_Q\nth(t)$ converges in distribution to $\M_Q(t)$,  $c_{P,Q}(\tilde{\M}_Q\nth(t))$ converges in distribution to $c_{P,Q}({\M}_Q(t))$, by the continuous mapping theorem. Combining this with \ref{lemma:CR.R1} and using Cram\'er-Slutzky's theorem proves \ref{lemma:CR.ntoC}.

Finally, \ref{lemma:CR.UI} holds
since the sequence 
$
c_{P,Q}(\tilde{\M}_Q\nth(t)) r\nth_{P,Q}(\tilde{\Y}\nth_Q(t), \y\nth_0)
$
converges in distribution to $c_{P,Q}({\M}_Q(t))$, as $n\to\infty$, by \ref{lemma:CR.ntoC},
and the  expectations are constantly equal to $1$, by \ref{lemma:CR.exp.n.1},\ref{lemma:CR.exp.1},  
thus by Remark  \ref{remark:thm}
the sequence is uniformly integrable.
\end{proof}
\end{lemma}

We now prove Theorem \ref{thm:weak_conv} for a general mutation matrix $P$, using properties \ref{lemma:CR.R1} and \ref{lemma:CR.UI} of Lemma \ref{lemma:CR}, Theorem \ref{thm:change_measure} and the PIM version of Theorem \ref{thm:weak_conv} (applied to  $Q$).

Fix a bounded continuous function $g:D_{\mathbb{R}_+^d \times \mathbb{N}^{d^2}} [0,t]\to \mathbb{R}$. To shorten the notation, let  
    $
    G\nth(t)=g(\{\tilde{\Z}_Q\nth(s)\}_{ 0\leq s\leq t})
    $
and 
    $
    G(t)=
    g(\{{\Z}_Q(s)\}_{ 0\leq s\leq t})
    $.
Because of Theorem \ref{thm:change_measure}, proving Theorem \ref{thm:weak_conv}  is equivalent to proving that 
\begin{equation}
\label{eq:convexpGRC}
\E{
    G\nth(t)
    r_{P,Q}\nth(\tilde{\Y}\nth_Q(t), \y\nth_0)
    c_{P,Q}(\tilde{\M}_Q\nth(t))
}
\xrightarrow[n\to\infty]{}
\E{
    G(t)
    c_{P,Q}({\M}_Q(t))
}.
\end{equation}
By the PIM version of  Theorem \ref{thm:weak_conv}  $\{\tilde{\Z}\nth_Q(s)\}_{0\leq s\leq t} \to \{ \Z_Q(s)\}_{0\leq s\leq t}$  in $D_{\mathbb{R}_+^d \times \mathbb{N}^{d^2}} [0,t]$.
Thus, since $g \cdot c_{P,Q}$ is a continuous function, 
the  continuous mapping theorem implies that
$
G\nth(t)c_{P,Q}(\tilde{\M}_Q\nth(t))
\to
G(t)c_{P,Q}({\M}_Q(t))
$, in distribution as $n\to\infty$.
Therefore, since $r\nth_{P,Q}(\tilde{\Y}\nth_Q(t), \y\nth_0) \to 1$ in probability,  as $n\to\infty$, by Lemma \ref{lemma:CR} \ref{lemma:CR.R1}, 
Cram\'{e}r- Slutzky's theorem implies
    \begin{align*}
    G\nth(t)
    c_{P,Q}(\tilde{\M}_Q\nth(t))
    r\nth_{P,Q}(\tilde{\Y}\nth_Q(t), \y\nth_0)
    \xrightarrow[n\to\infty]{}
    G(t)
    c_{P,Q}({\M}_Q(t)),
    \end{align*}
in distribution. 
Furthermore, 	
the sequence 
    $
    \{ G\nth(t)
    r\nth_{P,Q}(\tilde{\Y}\nth_Q(t), \y\nth_0)
    c_{P,Q}(\tilde{\M}_Q\nth(t))
    \}_{n\in\mathbb{N}}
    $
is uniformly integrable, because 
    $
    \{ 
    r\nth_{P,Q}(\tilde{\Y}\nth_Q(t), \y\nth_0)
    c_{P,Q}(\tilde{\M}_Q\nth(t))
    \}_{n\in\mathbb{N}}
    $
is uniformly integrable, by Lemma \ref{lemma:CR} \ref{lemma:CR.UI}, and $G\nth(t)$ is bounded.

Uniform integrability together with the convergence in distribution implies  that\\
$\E{|G\nth(t)| r\nth_{P,Q}(\tilde{\Y}\nth_Q(t), \y\nth_0) c_{P,Q}(\tilde{\M}_Q\nth(t))}$ 
converges to
$\E{|G(t)| c_{P,Q}({\M}_Q(t))}$ 
by Remark \ref{remark:thm}.
Therefore, if  $g$ is nonnegative Theorem \ref{thm:weak_conv} is proved, otherwise write $g$ as the difference between its positive and negative part and apply the previous argument on both parts.

The proof of Theorem \ref{thm:weak_conv} is completed. Note also that \eqref{eq:extended_weak_conv} is implied. In fact, with a similar argument to the one in the previous subsection, 
convergence 
of 
$\xi\tri_n(\{\tilde{\Z}_P\nth(s)\}_{0\leq s\leq t})$  to $\{\Z_P\tri(s)\}_{0\leq s\leq t} $
in $D_{E\tri}[0,t]$, 
is implied, since $t$ is fixed and $t<\norm{\y_0}$.
Furthermore, 
the convergence in $D_{E\tri}[0,t_1]$	for all $0<t_1<|\y_0|$ implies the
convergence in $D_{E\tri}[0,t_2)$ for all $t_2>|\y_0|$.
This holds by  e.g \cite[Ch. IV]{pollard1984},
defining the approximation map
$a:D_{E\tri}[0,t_2]\to D_{E\tri}[0,t_1]$, as
$a_{t_1} (z)(s)= z(s\wedge t_1)$ and observing that the previous part of the proof implies that, 
for all $ 0<t_1<|\y_0|$, $ a_{t_1}(\xi_n\tri(\tilde{\Z}_P\nth))\to a_{t_1}(\Z\tri_P) $ in $D_{E\tri}[0,t_2]$.
Therefore $\xi_n\tri(\tilde{\Z}\nth_P)$ converges to $ \Z\tri_P$ in $D_{E\tri}[0,t]$ for all $t>0$.  Convergence 	in $D_{E\tri}[0,t]$ for all $t>0$ implies convergence in $D_{E\tri}[0,\infty)$ by using for example \cite[Thm 16.7]{billingsley1999},  with the set of continuity points being $\mathbb{R}_{\geq0}$ since for all $t>0$  the probability that the limiting process $\Z\tri_P$ is discontinuous at $t$ is zero.

\begin{acks}
We would like to thank the  anonymous reviewers for  valuable comments which led to an improvement of the manuscript. 
MF is supported by the Knut and Alice Wallenberg Foundation (Program for Mathematics, grant 2020.072), HH is supported by the Swedish Research Council and MedTechLabs.

\end{acks}

\appendix
\section{Appendix}

\subsection{Backwards transition probabilities }
\label{appendix:transition_prob}


The backward transition probabilities of the block-counting jump chain of the coalescent have the form
    \begin{equation*}
    \label{eq:rho}
    \begin{aligned}
    \rho\nth (\vv\mid \y)
    =&
    \begin{cases}
    \frac{n y_i(n y_i -1)}{n\norm{\y}(n\norm{\y}-1+\theta)} \frac{1}{\pi(i\mid n\y -\e_i)}
    &\text{  if  }
    \vv=\e_i,
    \quad i=1\dots d,
    \\
    \frac{n y_i  \theta P_{ji}}{n\norm{\y}(n\norm{\y}-1+\theta)}
    \frac{\pi(j\mid n\y -\e_i)}{\pi(i\mid n\y -\e_i)}
    &\text{  if  }
    \vv=\e_i-\e_j,
    \quad i, j=1\dots d,
    \\
    0 &\text{  otherwise},
    \end{cases}
    \end{aligned}
    \end{equation*}
where $\pi(i\mid n\y)= \frac{n y_i +1}{n\norm{\y}+1}\frac{p(n\y+\e_i)}{p(n\y)}$ is 
the probability of sampling type $i$ after having sampled configuration $n\y$, and $p$ is the sampling probability described in Section \ref{sect:change_measure},
see e.g. \cite{deiorio2004}.
The expression above can be obtained from the expression of the forward transition probabilities by  employing Bayes' rule. 
In the PIM case, it is known that $\pi$ can be explicitly written as
    $$
    \pi(i\mid n\y)=\frac{n y_i +\theta Q_i}{n\norm{\y}+\theta },
    $$
which gives an explicit expression also for the backwards transition probabilities, given in \eqref{eq:rho_PIM}.

\subsection{Metric state-spaces and spaces of continuous functions} 
\label{appendix:spaces}
	
The goal of this section is to analyse the metric space $(E_1,\metric_1)$, in order to show that compact sets 
are bounded away from the boundary $\Omega_0$,  Remark \ref{remark:compact}, and to characterise continuous functions on the metric space $(E,\metric)$, Remark \ref{remark:continuous}, with the properties of  compact support and vanishing at infinity, Remark \ref{remark:vanishing_compact}. 
Characterisations in terms of the Euclidean metric are provided.

Let $S_1:= \mathbb{R}_+^d $ and $S:=S_1 \times \mathbb{N}^{d^2}$ be equipped with the Euclidean metric $\norma{\cdot}$.
Recall that $\Omega_0=\{\y=(y_1,\dots,y_d) : y_i=0 \text{ for some } i \}$ and  $\Omega_\infty:=\{\y=(y_1,\dots,y_d) : y_i=\infty \text{ for some } i \}$.
Let $\omega: E_1\to S_1$ be the homeomorphism defined by $\omega(\y)=\left(\frac{1}{y_1},\dots, \frac{1}{y_d}\right)$, with the convention $\frac{1}{\infty}=0$ and $\frac{1}{0}=\infty$. Recall that 
for  $\y_1, \y_2 \in E_1 $, 
    \begin{equation*}
    \metric_1(\y_1,\y_2)=
    \norma{\omega(\y_1)- \omega(\y_2) }.
    \end{equation*}
By equipping $E_1$ with the metric $\metric_1$, the roles of $0$ and $\infty$ are switched, 
whereas,  in the rest of the space, i.e. in  $ E_1 \setminus \Omega_\infty= S_1\setminus \Omega_0= (0,\infty)^d $,  the metric $\metric_1$ is equivalent to the Euclidean metric.
In fact,
if $\{\y\nth\}_{n=1}^{\infty}\subset (0,\infty)^d$ and $\x\in (0,\infty)^d$ then
$\metric_1(\y\nth,\x)\to 0 $ if and only if $\norma{\y\nth-\x}\to 0$.
This allows to characterise closed sets and compact sets  of $(E_1,\metric_1)$ in terms of the Euclidean metric.
	
\begin{remark}[Closed sets]
A set  $C\subset(E_1,\metric_1)$ is closed  if and only if 
$\omega(C)=\{\textbf{x}\in S_1: \textbf{x}=\frac{1}{\y} \text{ for some } \y \in E_1\}$ is closed in $(S_1,\norma{\cdot})$. Equivalently, $ C$ is closed  if and only if one of the following conditions is satisfied
\begin{itemize}
\item
$C\cap \Omega_\infty=\emptyset$ and 
$C$ is bounded and closed in $(E_1\setminus\Omega_\infty , \norma{\cdot})$ ;
\item
$C\setminus \Omega_\infty$ is closed in $(E_1\setminus\{\infty\} , \norma{\cdot})$  and $\y\in C\cap \Omega_\infty$ for all $\y\in \Omega_\infty$ such that there exists  a sequence $\{\y\nth\}_{n\in \mathbb{N}}\subset C$, with $y\nth_i\to y_i=\infty$, for  $i\in \mathcal{I}\neq \emptyset$, and $y\nth_j\to y_j\in (0,\infty) $, for $j\in\{1,\dots,d\}\setminus \mathcal{I} $.
\end{itemize}
\end{remark}

\begin{remark}[Compact sets]
\label{remark:compact}
A compact set in $(E_1,\metric_1)$ is a closed set which is bounded away from $\Omega_0$. More precisely, a set $K\subset (E_1,\metric_1)$ is compact if and only if
$K$ is closed in $(E_1,\metric_1)$ and there exists $\delta>0$ such that, for all  $ \x\in K$,
$y_j>\delta,$  $j=1,\dots,d $.
\end{remark}	
Furthermore,  real-valued continuous functions on $(E_1,\metric_1)$ are  simply continuous functions on $(E_1\setminus\Omega_\infty,\norma{\cdot})$  that can be continuously extended on  $\Omega_\infty$, 
as specified in the following.	
\begin{remark}[Continuous functions]
\label{remark:continuous}
A function $f: (E_1,\metric_1)\to \mathbb{R}$ is continuous if and only if the two following conditions are satisfied
\begin{itemize}
\item 
its restriction $f:(E_1\setminus \Omega_\infty,\norma{\cdot}) \to \mathbb{R}$ is continuous   
\item
$f(\y\nth)\to f(\y)$ for any sequence
$\{\y\nth\}_{n\in \mathbb{N}}\subset E_1\setminus \Omega_\infty$, with $y\nth_i\to y_i=\infty$, for $i\in \mathcal{I}\neq \emptyset$, and $y\nth_j\to y_j\in (0,\infty) $, for $j\in\{1,\dots,d\}\setminus \mathcal{I} $.
\end{itemize}
\end{remark}
Since $E=E_1 \times \mathbb{N}^{d^2}$ is equipped with the product metric $\metric= \metric_1 \oplus \left\lVert \cdot\right\rVert_2 $, and $\mathbb{N}^{d^2}$ is a discrete space, the above characterisations can be easily extended to $(E,\metric)$. We focus now on spaces of continuous functions on $(E,\metric)$.
In particular, a continuous function on $(E,\metric)$ is simply a function which is continuous with respect to its component in $E_1$ in the sense of Remark \ref{remark:continuous}.
We say that  a function  vanishes at infinity if, 
given any  $\epsilon > 0 $, 
there exists a compact set $K_{\epsilon}\subset E$ such that  $\left|f(\z)\right|<\epsilon$, $\forall \z\in K_\epsilon^c$.
Since the component in $(E_1,\metric_1)$ of a 
compact set is bounded away from $\Omega_0$ and the component in   $\mathbb{N}^{d^2}$ is bounded,
the following  holds. 
\begin{remark}[Vanishing at infinity, compact support]
\label{remark:vanishing_compact}
A function $f:(E,\metric)\to \mathbb{R}$ 
\begin{itemize}
\item
vanishes at infinity if and only if, for any given $\epsilon>0$, 
there exist $\delta_{\epsilon},M_{\epsilon}>0$ such that 
    \begin{equation*}
    \begin{cases}
    y_j<\delta_{\epsilon} \text{ for  some } j=1,\dots,d \\
    \text{or} \\
    m_{ij}>M_{\epsilon} \text{ for  some }  i,j=1,\dots,d
    \end{cases}
    \Rightarrow
    |f(\x,\m)|<\epsilon
    \end{equation*}
\item
has compact support if and only if there exist $\delta,M>0$ such that 
    \begin{equation*}
    \begin{cases}
    y_j<\delta \text{ for  some } j=1,\dots,d \\
    \text{or} \\
    m_{ij}>M \text{ for  some }  i,j=1,\dots,d
    \end{cases}
    \Rightarrow
    f(\x,\m)=0
    \end{equation*}
\end{itemize}
\end{remark}	
The characterisations of this section provide a description of the spaces of continuous functions on $(E,\metric)$ vanishing at infinity, $\hat{C}(E)$, and of smooth functions on $(E,\metric)$ with compact support, $C_c^{\infty}(E)$. 
Furthermore, note that $C_c^{\infty}(E)$ is dense in $\hat{C}(E)$ which is dense in the  space of bounded continuous functions on $(E,\metric)$. Moreover, all functions in $C_c^{\infty}(E)$ and $\hat{C}(E)$ are uniformly continuous, while functions in $C_c^{\infty}(E)$, and all of their derivatives, are  Lipschitz continuous.

Finally, we characterise continuous functions on the extended state-space $E^\Delta = E \cup \{\Delta\}$, which is the one-point,  or Alexandroff, compactification of $E$ equipped with the  standard  topology inherited by the topology in $E$.
Since $E $ is a separable locally compact metric space, its one-point compactification is metrizable, see e.g. \cite{mandelkern1989}. Let $(E\tri,\metric\tri)$ be the corresponding metric space. 

\begin{remark}[Extended continuous functions]
\label{remark:continuous_extended}
A function $f: (E\tri,\metric\tri)\to \mathbb{R}$ is continuous if and only if the  following conditions are satisfied
\begin{itemize}
\item 
its restriction $f:((E_1\setminus \Omega_\infty)\times \mathbb{N}^{d^2},\norma{\cdot}) \to \mathbb{R}$ is continuous;  
\item
$\forall \m \in \mathbb{N}^{d^2}, f(\y\nth,\m)\to f(\y,\m)$ for any sequence 
$\{\y\nth\}_{n\in \mathbb{N}}\subset E_1\setminus \Omega_\infty$, with $y\nth_i\to y_i=\infty$, for  $i\in \mathcal{I}\neq \emptyset$, and $y\nth_j\to y_j\in (0,\infty) $, for $j\in\{1,\dots,d\}\setminus \mathcal{I} $;
\item
$f(\y\nth,\m\nth)\to f(\Delta)$ for any sequence $\{(\y\nth,\m\nth\}_{n\in \mathbb{N}}\subset (E_1\setminus \Omega_\infty)\times \mathbb{N}^{d^2} $, with $y\nth_i\to 0$, for some $i\in\{1,\dots,d\}$ or $m_{ij}\nth\to \infty$ for some $i,j\in\{1,\dots,d\}$.
\end{itemize}
\end{remark}
Note that, while using the metrics $\psi$ and $\psi\tri$ allows for a more compact formulation,
 all the statements concerning continuous functions on $E$ or $E\tri$, can be expressed in terms of the Euclidean metric by using the remarks of this section.

\subsection{Infinitesimal generator and semigroup of the limiting process (PIM)}
\label{appendix:generator_semigroup}


In this section, under the assumption of parent independent mutations, we show that the semigroup $T$, defined in \eqref{eq:Tn}, is a positive, strongly continuous, contraction semigroup and that the infinitesimal generator $A$, associated to the semigroup $T$, is indeed of the form \eqref{eq:A}.
The following, straightforward to derive,  expressions for the probabilities in \eqref{eq:gamma_prob} are needed
    \begin{equation}
    \label{p0}
    \gamma_{\zero}(t,\x)=\left(\frac{\norm{\x}-t}{\norm{\x}}\right)^{\theta} ,
    \end{equation}	
    \begin{equation}
    \label{pij}
    \gamma_{\e_{ij}}(t,\x)=\left(\frac{\norm{\x}-t}{\norm{\x}}\right)^{\theta} 
    \frac{ \theta Q_j y_i }{\norm{\x}}
    \log\left(\frac{\norm{\x}}{\norm{\x}-t}\right),
    \quad i,j=1,\dots,d.
    \end{equation}	
Furthermore,
since  the sum of independent Poisson random variables is Poisson distributed,  for all  	$k\in \mathbb{N}$,
    \begin{equation}
    \label{sumpw}
    \sum_{\w \in \mathbb{N}^{d^2}: \norm{\w}\geq k} \gamma_{\w}(t,\x)=
    \Prob{\norm{\M(t+s)-\M(s)}\geq k}
    =
    \sum_{h=k}^{\infty} \frac{\Lambda(t,\x)^h}{h!} e^{-\Lambda(t,\x)},
    \end{equation}
where 
    \begin{equation*}
    \Lambda(t,\x)=\sum_{i,j=1}^{d} \Lambda_{ij}(t,\x)= 
    \theta \log\left(\frac{\norm{\x}}{\norm{\x}-t}\right).
    \end{equation*}

\subsubsection*{Semigroup}
$T$ is a positive  contraction semigroup, 
being  the semigroup associated to a Markov process.
It remains to show that, for all $f\in \hat{C}(E)$,
    \begin{equation*}
    \lim_{t\to0}	\sup_{(\x,\m)\in E}
    \left| 
    T(t)f(\x,\m)-f(\x,\m)	
    \right|
    =0.
    \end{equation*} 
Take  $f\in \hat{C}(E)$, $f$ is continuous and  vanishing in the sense of Remark \ref{remark:vanishing_compact},  thus it is  uniformly continuous with respect to the metric $\metric$  and bounded. Given $\epsilon>0$, let $K_\epsilon=\{(\x,\m)\in E: y_j\geq \delta_\epsilon, m_{ij}\leq M_\epsilon, \forall i,j=1,\dots,d\}$,
for some $\delta_\epsilon,M_\epsilon>0$, be the compact set such that 
$|f(\z)|<\epsilon$, for all $\z\in K_\epsilon^c$. 
The supremum in the last display is analysed in the three different  cases;
$\x\in E_1\cap \Omega_\infty$, and, when $\x\in E_1\setminus \Omega_\infty$, $\norm{\x} \leq t$ and $\norm{\x} > t$.

If $\x\in E_1\cap \Omega_\infty$, then 
$T(t)f(\x,\m)-f(\x,\m)	=0,  \forall \m \in  \mathbb{N}^{d^2} ,\forall t\geq 0$.
Assume now
$\x\in E_1\setminus \Omega_\infty$. 
If $\norm{\x}\leq t$, then
$T(t)f(\x,\m)=0$, 
and 
for $t<\delta_\epsilon$, 
$\norm{\x}\leq t<\delta_\epsilon$ implies $(\x,\m)\in K_\epsilon^c$ and thus $|f(\x,\m)|<\epsilon$. Therefore, 
\begin{align*}
\sup_{\substack{(\x,\m)\in E\\\norm{\x}\leq t \\ \x\notin \Omega_\infty}}
\left| 
T(t)f(\x,\m)
-
f(\x,\m)	
\right|	
< \epsilon .
\end{align*}
If $\norm{\x}> t$, then,  by \eqref{semigroupT},
\begin{align*}
|
T(t)f(\x,\m)	
-
f(\x,\m)|=
\left| 
\sum_{\w \in \mathbb{N}^{d^2}}
f\left(\x- \frac{\x}{\norm{\x}}t,\m+\w\right) \gamma_{\w}(t,\x)
-	
f(\x,\m)	
\right|.
\end{align*}
The remaining argument is divided into the cases $(\x,\m)\in K_\epsilon^c$ and $(\x,\m)\in K_\epsilon$. If $(\x,\m)\in K_\epsilon^c$, then $\left(\x- \frac{\x}{\norm{\x}}t,\m+\w\right)\in K_\epsilon^c$ and 
$|f(\x,\m)|, f\left(\x- \frac{\x}{\norm{\x}}t,\m+\w\right)<\epsilon $.
Since 
$\sum_{\w \in \mathbb{N}^{d^2}} \gamma_{\w}(t)=1$, it follows that
\begin{align*}
\sup_{\substack{(\x,\m)\in K_\epsilon^c \\ \norm{\x}> t \\ \x\notin \Omega_\infty }}
\left| 
T(t)f(\x,\m)
-
f(\x,\m)	
\right|	
< 2 \epsilon.
\end{align*}
If $(\x,\m)\in K_\epsilon$,  consider the following inequality,  
\begin{align*}
&
\left| 
\sum_{\w \in \mathbb{N}^{d^2}} 
f\left(\x- \frac{\x}{\norm{\x}}t,\m+\w\right) \gamma_{\w}(t,\x)
-	
f(\x,\m)\right| 
\\ 
&
\quad\leq
\left| 
f\left(\x- \frac{\x}{\norm{\x}}t,\m\right)
-	
f(\x,\m)	
\right|
\gamma_{\zero}(t,\x)
\\
& \qquad +
\left|  f(\x,\m)  \right|
\left( 1-\gamma_{\zero}(t,\x) \right)
\\
& \qquad+
\sum_{\substack{\w \in \mathbb{N}^{d^2}\\ \norm{\w}\geq 1}}
\left| f\left(\x- \frac{\x}{\norm{\x}}t,\m+\w\right) \right|  \gamma_{\w}(t,\x) .
\end{align*}
Since $f$ is uniformly continuous,   there exists $\eta_\epsilon>0$, not depending on $\x$, such that \\
$
\left| 	f\left(\x- \frac{\x}{\norm{\x}}t,\m\right) -	f(\x,\m)	\right|<\epsilon,
$ 
if 
$
\metric\left(
\left(\x- \frac{\x}{\norm{\x}}t,\m\right),(\x,\m)
\right) <\eta_\epsilon
$.
Keeping in mind that  $\norm{\x}$,$\norma{\x}>\delta_\epsilon$, note that
\begin{align}
\label{bounddistance}
\metric\left(
\left(\x- \frac{\x}{\norm{\x}}t,\m\right),(\x,\m)
\right)
= 	
\metric_1  \left(\x- \frac{\x}{\norm{\x}}t,\x\right) 
= 
\frac{t}{\norm{\x}-t} \sum_{i=1}^d\frac{1}{y_i}
<
\frac{t}{\delta_\epsilon -t} \frac{d}{\delta_\epsilon}
\end{align}
is bounded by $\eta_\epsilon$,	
if $t< \eta_\epsilon '  $, where $\eta_\epsilon '>0$ only depends on $\delta_\epsilon,\eta_\epsilon$, not  on $\x$.
Recalling that $f$ is bounded  and
using (\ref{p0}) yields
\begin{equation}
\label{bound1-p0}
\left|  f(\x,\m)  \right|
\left (1-\gamma_{\zero}(t,\x) \right) 
\leq \norminfty{f}
\left( 1-  \left( 1- \frac{t}{\delta_\epsilon}\right)^\theta \right), 
\end{equation}
which is bounded by $\epsilon$, 
if $t < \eta _\epsilon '' $, for some  $\eta_\epsilon ''>0$ which only depends on $\delta_\epsilon$, not  on $\x$.
Furthermore, using (\ref{sumpw}) and again that $f$ is bounded  yields
\begin{align}
\allowdisplaybreaks
\label{boundsumpw}
\nonumber
\sum_{\substack{\w \in \mathbb{N}^{d^2}\\ \norm{\w}\geq 1}} &
\left| f\left(\x- \frac{\x}{\norm{\x}}t,\m+\w\right) \right|  \gamma_{\w}(t,\x) 
\\ \nonumber & \leq
\norminfty{f}
\sum_{\substack{\w \in \mathbb{N}^{d^2}\\ \norm{\w}\geq 1}}
\gamma_{\w}(t,\x) 
\\\nonumber &=
\norminfty{f}
\left(\frac{\norm{\x}-t}{\norm{\x}}\right)^{\theta} 
\sum_{h=1}^{\infty}
\left[   \theta \log\left(\frac{\norm{\x}}{\norm{\x}-t}\right) \right]^h \frac{1}{h!}	 
\\
&\leq
\norminfty{f}
\sum_{h=1}^{\infty}
\left(    \frac{\theta t}{\delta_\epsilon-t} \right)^h \frac{1}{h!},	
\end{align}
where the last quantity is bounded by $\epsilon$ if $t<\eta_\epsilon '''$,  for some $\eta_\epsilon '''>0$, which only depends on $\delta_\epsilon$, not  on $\x$.
Therefore, 
if $t< \min\{ \eta_\epsilon', \eta_\epsilon'', \eta_\epsilon''' \}$,
\begin{align*}
\sup_{\substack{(\x,\m)\in E\\\norm{\x}>t \\ \x\notin \Omega_\infty}}
\left| 
T(t)f(\x,\m)
-
f(\x,\m)	
\right|	
< 3 \epsilon,
\end{align*}
which completes the proof.

\subsubsection*{Infinitesimal generator}
Let $A$ be the operator defined as in \eqref{eq:A}, for $f\in C^\infty_c(E)$, 
	\begin{equation*}
	A f(\x,\m)= 
	-\left\langle 
	\nabla_\x f(\x,\m),\frac{\x}{\norm{\x}}
	\right\rangle
	+
	\sum_{i,j=1}^d
	\left[
	f(\x,\m+\e_{ij})-f(\x,\m)
	\right]
	\lambda_{ij}(\y),
	\end{equation*}
We show that $A$ is the infinitesimal generator of the process  $\Z$, associated to the semigroup $T$, by proving that
    \begin{align*}
    \lim_{t\to0}	\sup_{(\x,\m)\in E}
    \left| 
    \frac{T(t)f(\x,\m)-f(\x,\m)}{t}	-Af(\x,\m)
    \right|
    =0 .
    \end{align*}

Take $f\in C^\infty_c(E)$,  $f $ is smooth  with 
compact support in the sense of Remark \ref{remark:vanishing_compact}. Thus $f$ and all of its derivatives are  Lipschitz continuous with respect to $\metric$ and bounded.   
Let  $K=\{(\x,\m)\in E : y_j\geq \delta, m_{ij}\leq M, \forall i,j=1,\dots,d\}$ 
for some $\delta,M>0$, be the compact set that contains the support of $f$, so that  $f(\x,\m)=0$,
$\forall (\x,\m)\in K^c$.

Since the limit as $t\to0$ is considered, assume $t<\delta$.
Three cases are possible: $\x\in E_1\cap \Omega_\infty$, and when $\x\in E_1\setminus \Omega_\infty$, $\norm{\x} \leq t$ and $\norm{\x} > t$.
If $\x\in E_1\cap \Omega_\infty$, then $T(t)f(\x,\m)=f(\x,\m)$ and $Af(\x,\m)=0$.
Assume now $\x\in E_1\setminus \Omega_\infty$.
If $\norm{\x}\leq t$,  then $T(t)f(\x,\m)=0$, and  
$\norm{\x}\leq t<\delta$ implies that $(\x,\m) \in K^c$,  $(\x,\m+\e_{ij}) \in K^c, i,j=1,\dots,d$, and that $f$ is equal to $0$ in a neighbourhood of $(\x,\m)$, hence $Af(\x,\m)=0$. 
If $\norm{\x} > t$ and $(\x,\m)\in K^c$, then, as above, $f(\x,\m)=0$ and $Af(\x,\m)=0$. 
Furthermore $\left(\x- \frac{\x}{\norm{\x}}t,\m+\w\right)\in K^c, \forall \w \in \mathbb{N}^{d^2}$, 
and thus
$
T(t)f(\x,\m)=0
$

Consider now  the supremum  over all $(\x,\m)\in K$, with $ \norm{\x}>t$, $\x\in E_1\setminus \Omega_\infty$.
The aim is to bound the absolute value in the new display (\ref{boundsA})   by some function depending on $t$ and $\delta$, not on $\x$ and $\m$, which vanishes as $t\to 0$. Start by considering the bound,
\begin{align}
\allowdisplaybreaks
\nonumber
&
\left| 
\frac{T(t)f(\x,\m)-  f(\x,\m)}{t}	-Af(\x,\m)
\right| 
\\ \nonumber
&\leq
\left| 
\frac{f\left(\x- \frac{\x}{\norm{\x}}t,\m\right)-f(\x,\m)}{t}	
\gamma_{\zero}(t,\x)
+
\left\langle 
\nabla_\x f(\x,\m),\frac{\x}{\norm{\x}}
\right\rangle
\right| 
\\ \nonumber
&+
\sum_{i,j=1}^d
\left|
\frac{f\left(\x- \frac{\x}{\norm{\x}}t,\m+\e_{ij}\right)-f(\x,\m)}{t}	
\gamma_{\e_{ij}}(t,\x)
-
\left[
f(\x,\m+\e_{ij})-f(\x,\m)
\right]
\lambda_{Q,ij}(\y)
\right| 
\\	\label{boundsA}
&+
\sum_{\w\in\mathbb{N}^{d^2},\norm{w}\geq 2}
\left|
\frac{f\left(\x- \frac{\x}{\norm{\x}}t,\m+\w\right)-f(\x,\m)}{t}	
\right| 
\gamma_{\w}(t,\x).
\end{align}
Each of the three terms in the right hand side above are studied separately.

To bound the first term, apply the mean value theorem,  
\begin{align*}
f\left(\x- \frac{\x}{\norm{\x}}t,\m\right)-f(\x,\m)= 
\left\langle 
\nabla_\x f\left(a\x +(1-a)\left(\x- \frac{\x}{\norm{\x}}t \right),\m\right),-\frac{\x}{\norm{\x}}t
\right\rangle
\end{align*}
for some $a\in(0,1)$, to obtain
\begin{align}
\allowdisplaybreaks
\label{boundA1} \nonumber
&
\left| 
\frac{f\left(\x- \frac{\x}{\norm{\x}}t,\m\right)-f(\x,\m)}{t}	
\gamma_{\zero}(t,\x)
+
\left\langle 
\nabla_\x f(\x,\m),\frac{\x}{\norm{\x}}
\right\rangle
\right| 	
\\ \nonumber
&\qquad \leq 
\left| 
\left\langle 
\nabla_\x f\left(( \norm{\x}-(1-a)t)\frac{\x}{\norm{\x}} ,\m\right)
-\nabla_\x f (\x,\m)
,\frac{\x}{\norm{\x}}
\right\rangle
\right| 	   \gamma_{\zero}(t,\x)	 \\
&\qquad +
\left| \left\langle 
\nabla_\x f\left(\x ,\m\right)
,\frac{\x}{\norm{\x}}
\right\rangle  \right| 	 
(1- \gamma_{\zero}(t,\x)	) .
\end{align}
Note that  $\nabla_\x f$ is bounded and Lipschitz continuous  and
\begin{align*}
\metric_1\left(( \norm{\x}-(1-a)t)\frac{\x}{\norm{\x}} ,\x\right)
=
\frac{(1-a)t}{\norm{\x}-(1-a)t} \sum_{i=1}^d \frac{1}{y_i}
\leq
\frac{t}{\delta-(1-a)t}\frac{d}{ \delta}.
\end{align*}
Using  that  $\frac{\x}{\norm{\x}}\leq 1$ and $\norma{\x},\norm{\x}>\delta$, and   (\ref{p0}), (\ref{bound1-p0}) for $\gamma_\zero  $ and $1-\gamma_\zero$, 
it follows that
(\ref{boundA1}) is bounded by
\begin{equation*}
\begin{aligned}
\norminfty{\nabla^2_\x f}
\frac{d t}{\delta(\delta-(1-a)t)}
+
\norminfty{\nabla_\x f}
\left(
1-  \left( 1- \frac{t}{\delta}\right)^\theta \right).
\end{aligned}
\end{equation*}
Next, consider the $ij^{th}$ term of the sum in  (\ref{boundsA}), which can be bounded by, 
\begin{align}
\label{boundA2}
\allowdisplaybreaks \nonumber
&
\left|
\frac{f\left(\x- \frac{\x}{\norm{\x}}t,\m+\e_{ij}\right)-f(\x,\m)}{t}	
\gamma_{\e_{ij}}(t,\x)
-
\left[
f(\x,\m+\e_{ij})-f(\x,\m)
\right]
\frac{\theta Q_j y_i}{\norm{\x}^2}
\right|  \\ \nonumber
&\quad\leq
\left|
f\left(\x- \frac{\x}{\norm{\x}}t,\m+\e_{ij}\right)-f(\x,\m+\e_{ij})
\right|
\frac{1}{t}
\gamma_{\e_{ij}}(t,\x)
\\
&\qquad +
\left| f(\x,\m+\e_{ij})-f(\x,\m) \right|
\left|\frac{ \gamma_{\e_{ij}}(t,\x)}{t} - \frac{\theta Q_j y_i}{\norm{\x}^2} \right| .
\end{align}
Recalling that
$f$ is Lipschitz continuous, using the bound  (\ref{bounddistance}) for the distance  and  expression (\ref{pij}) for $\gamma_{\e_{ij}}$, yields, 
\begin{align*}
\allowdisplaybreaks
&
\frac{1}{t}
\gamma_{\e_{ij}}(t,\x)
\left|
f\left(\x- \frac{\x}{\norm{\x}}t,\m+\e_{ij}\right)-f(\x,\m+\e_{ij})
\right|	
\\
&\quad\leq
-\frac{1}{t}\theta Q_j  \log \left( 1-\frac{t}{\delta}  \right)
\norminfty{\nabla_\x f} \ \metric_1\left(\x- \frac{\x}{\norm{\x}}t,\x\right)  \\
&\quad<
- \theta\norminfty{\nabla_\x f}    \frac{1}{t}\log \left( 1-\frac{t}{\delta}  \right)
\frac{t d}{\delta(\delta -t)} .
\end{align*}
Using the explicit expressions for the probabilities  (\ref{p0}), (\ref{pij}), and  the bound (\ref{bound1-p0}) for $1-\gamma_\zero$, yields, 
\begin{align*}
\allowdisplaybreaks
&
\left| f(\x,\m+\e_{ij})-f(\x,\m) \right|
\left| \gamma_{\e_{ij}}(t,\x) - \frac{\theta Q_j y_i}{\norm{\x}^2} \right|   \\
&\quad\leq
2 \norminfty{ f}
\frac{\theta Q_j y_i}{\norm{\x}} 
\left(
\gamma_\zero(t,\x) 
\left|  \frac{1}{t}\log \left( \frac{\norm{\x}-t}{\norm{\x}}  \right)+\frac{1}{\norm{\x}}  \right|
+
\frac{1}{\norm{\x}}  (1-\gamma_\zero(t,\x))
\right)
\\
& \quad\leq
2 \theta  \norminfty{ f}
\left(
-\frac{1}{t}\log \left( \frac{\delta- t}{\delta }  \right)-\frac{1}{\delta}
+\frac{1}{\delta} \left( 1-  \left( 1- \frac{t}{\delta}\right)^\theta \right)
\right).
\end{align*}
Therefore, each term in (\ref{boundA2}) is bounded by a function of $t$ and $\delta$, not depending on $\x$ which vanishes as $t\to 0$.

Similarly to (\ref{boundsumpw}), 
the last sum in (\ref{boundsA})	is bounded as  follows
\begin{equation*}
\sum_{\w\in\mathbb{N}^{d^2},\norm{w}\geq 2}
\left|
\frac{f\left(\x- \frac{\x}{\norm{\x}}t,\m+\w\right)-f(\x,\m)}{t}	
\right| 
\gamma_{\w}(t, \x)
\leq
2 \norminfty{ f}
\frac{1}{t}
\sum_{h=2}^{\infty}
\left(    \frac{\theta t}{\delta-t} \right)^h \frac{1}{h!},	
\end{equation*}
which also vanishes as $t\to 0 $.	The  proof is now complete.



\begin{thebibliography}{10}

\bibitem{barbour2000}
A.~D. Barbour, S.~N. Ethier, and R.~C. Griffiths.
\newblock A transition function expansion for a diffusion model with selection.
\newblock {\em Annals of Applied Probability},
10(1):123--162, 2000.
 
\bibitem{bhaskar2014}
A.~{Bhaskar}, A.~G. {Clark}, and Y.~S. {Song}.
\newblock Distortion of genealogical properties when the sample is very large.
\newblock {\em Proceedings of the National Academy of Sciences of the United States of America}, 111(6):2385--2390, 2014.
	
\bibitem{billingsley1999}
P.~Billingsley.
\newblock {\em Convergence of probability measures}.
\newblock Wiley Series in Probability and Statistics. John Wiley \& Sons Inc.,
New York, second edition, 1999.

\bibitem{birkner2008}
M.~Birkner and J.~Blath.
\newblock Computing likelihoods for coalescents with multiple collisions in the
infinitely many sites model.
\newblock {\em Journal of mathematical biology}, 57(3):435--465, 2008.

\bibitem{birkner2011}
M.~Birkner, J.~Blath, and M.~Steinrücken.
\newblock Importance sampling for lambda-coalescents in the infinitely many
sites model.
\newblock {\em Theoretical population biology}, 79(4):155--173, 2011.

\bibitem{deiorio2004}
M.~{De~Iorio} and R.~C. Griffiths.
\newblock Importance sampling on coalescent histories. {I}.
\newblock {\em Advances in Applied Probability}, 36(2):417–433, 2004.





\bibitem{Etheridge2009}
A.~M. Etheridge and R.~C. Griffiths.
\newblock A coalescent dual process in a {M}oran model with genic selection.
\newblock {\em Theoretical Population Biology}, 75:320--330, 2009.

\bibitem{Ethier1986}
S.~N. Ethier and T.~G. Kurtz.
\newblock {\em Markov processes: characterization and convergence}, volume 282.
\newblock John Wiley \& Sons, 1986.

\bibitem{favero2021}
M.~Favero, H.~Hult, and T.~Koski.
\newblock A dual process for the coupled wright-fisher diffusion.
\newblock {\em Journal of Mathematical Biology}, 82(6), 2021.

\bibitem{favero2020}
M.~Favero and H.~Hult.
\newblock {Asymptotic behaviour of sampling and transition probabilities in
    coalescent models under selection and parent dependent mutations}.
\newblock {\em Electronic Communications in Probability}, 27:1 -- 13, 2022.

\bibitem{favero2023}
M.~Favero and P.~A.~Jenkins.
\newblock {Sampling probabilities, diffusions, ancestral graphs, and duality under strong selection}.
\newblock {\em Electronic Journal of Probability}, 30:1--43, 2025.



\bibitem{griffiths1991}
R.~C. Griffiths.
\newblock The two-locus ancestral graph.
\newblock In I.~V. Basawa and R.~L. Taylor, editors, {\em Selected Proceedings
    of the Sheffield Symposium on Applied Probability}, volume~18 of {\em Lecture
    Notes--Monograph Series}, pages 100--117. Institute of Mathematical
Statistics, Hayward, CA, 1991.

\bibitem{griffiths2008}
R.~C. Griffiths, P.~A. Jenkins, and Y.~S. Song.
\newblock Importance sampling and the two-locus model with subdivided
population structure.
\newblock {\em Advances in applied probability}, 40(2):473--500, 2008.

\bibitem{griffiths1997}
R.~C. Griffiths and P.~Marjoram.
\newblock An ancestral recombination graph.
\newblock In P.~Donelly and S.~Tavar\`e, editors, {\em Progress in Population
    Genetics and Human Evolution}, page 257–270. Springer-Verlag, Berlin, 1997.

\bibitem{griffiths1994simulating}
R.~C. {Griffiths} and S.~{Tavar\'e}.
\newblock Simulating probability distributions in the coalescent.
\newblock {\em Theoretical Population Biology}, 46(2):131--159, 1994.

\bibitem{jenkins2016}
R. C. Griffiths, P. A. Jenkins and S. Lessard.
\newblock A coalescent dual process for a Wright–Fisher diffusion with recombination and its application to haplotype partitioning.
\newblock {\em Theoretical Population Biology},
112: 126-138,
2016.

\bibitem{hobolth2008}
A.~Hobolth, M.~K. Uyenoyama, and C.~Wiuf.
\newblock Importance sampling for the infinite sites model.
\newblock {\em Statistical applications in genetics and molecular biology},
7(1):Article32, 2008.


\bibitem{kelleher2016}
J.~{Kelleher}, A.~M. {Etheridge}, and G.~{McVean}.
\newblock Efficient coalescent simulation and genealogical analysis for large
sample sizes.
\newblock {\em PLOS Computational Biology}, 12(5), 2016.

\bibitem{kingman1982b}
J.~F.~C. Kingman.
\newblock The coalescent.
\newblock {\em Stochastic Processes and their Applications}, 13(3):235 -- 248,
1982.



\bibitem{Koskela2015}
J.~Koskela, P.~A. Jenkins, and D.~Spanò.
\newblock Computational inference beyond {Kingman's} coalescent.
\newblock {\em Journal of Applied Probability}, 52(2):519--537, 06 2015.

\bibitem{Koskela2018}
J.~Koskela, D.~Span{\`o}, and P.~A. Jenkins.
\newblock Inference and rare event simulation for stopped {Markov processes via
    reverse-time sequential Monte Carlo}.
\newblock {\em Statistics and Computing}, 28(1):131--144, 2018.

\bibitem{Krone1997}
S.~M. Krone and C.~Neuhauser.
\newblock Ancestral processes with selection.
\newblock {\em Theoretical Population Biology}, 51:210--237, 1997.

\bibitem{kurtz2001}
T.~G. Kurtz.
\newblock {\em Lectures on Stochastic Analysis}.
\newblock University of Wisconsin - Madison, 2001.

\bibitem{mandelkern1989}
M.~Mandelkern.
\newblock Metrization of the one-point compactification.
\newblock {\em Proceedings of the AMS}, 107(4), 1989.

\bibitem{mohle2001}
M.~Möhle and S.~Sagitov.
\newblock A classification of coalescent processes for haploid exchangeable
population models.
\newblock {\em Annals of Probability}, 29:1547--1562, 2001.



\bibitem{Neuhauser1997}
C.~Neuhauser and S.~M. Krone.
\newblock The genealogy of samples in models with selection.
\newblock {\em Genetics}, 154:519--534, 1997.

\bibitem{pitman1999}
J.~Pitman.
\newblock Coalescent with multiple collisions.
\newblock {\em Annals of Probability}, 27:1870--1902, 1999.

\bibitem{pollard1984}
D.~Pollard.
\newblock {\em Convergence of Stochastic Processes}.
\newblock Springer Series in Statistics. Springer-Verlag New York, 1984.


\bibitem{sagitov1999}
S.~Sagitov.
\newblock The general coalescent with asynchronous mergers of ancestral lines.
\newblock {\em Journal of Applied Probability}, 36:1116--1125, 1999.

\bibitem{shiga1981}
T. Shiga.
\newblock Diffusion Processes in Population Genetics.
\newblock {\em J. Math. Kyoto Univ. (JMKYAZ)}, 21(1):133--151, 1981.

\bibitem{schweinsberg2000}
J.~Schweinsberg.
\newblock Coalescents with simultaneous multiple collisions.
\newblock {\em Electronic Journal of Probability}, 5:50pp, 2000.


\bibitem{stephens2000}
M.~Stephens and P.~Donnelly.
\newblock Inference in molecular population genetics.
\newblock {\em Journal of the Royal Statistical Society: Series B (Statistical
    Methodology)}, 62(4):605--635, 2000.

\bibitem{stephens2003}
M.~Stephens and P.~Donnelly.
\newblock Ancestral inference in population genetics models with selection
(with discussion).
\newblock {\em Australian \& New Zealand Journal of Statistics},
45(4):395--430, 12 2003.


\bibitem{wakeley2020}
J. Wakeley.
\newblock Developments in coalescent theory from single loci to chromosomes.
\newblock {\em Theoretical population biology},
133:56--64, 2020.



\end{thebibliography}
\end{document}